\documentclass[a4paper,11pt]{amsart}

\usepackage{fullpage,amsxtra,amssymb,amsthm,amsmath,amscd}
\usepackage[latin1]{inputenc}
\numberwithin{equation}{section}

\renewcommand{\leq}{\leqslant}
\renewcommand{\geq}{\geqslant}

\def\stacksum#1#2{{\stackrel{{\scriptstyle #1}}
{{\scriptstyle #2}}}}

\newcommand{\Cc}{\mathbf{C}}

\newcommand{\Aa}{\mathbf{A}}
\newcommand{\Zz}{\mathbf{Z}}

\newcommand{\Qq}{\mathbf{Q}}
\newcommand{\Fp}{\mathbf{F}}

\newcommand{\mods}[1]{\,(\mathrm{mod}\,{#1})}

\newcommand{\ideal}[1]{\mathfrak{{#1}}}

\newcommand{\ra}{\rightarrow}
\newcommand{\lra}{\longrightarrow}

\newcommand{\fleche}[1]{\stackrel{#1}{\lra}}

\newcommand{\barre}[1]{\overline{{#1}}}

\DeclareMathOperator{\spec}{Spec}

\DeclareMathOperator{\rank}{rank}

\DeclareMathOperator{\Reel}{Re}

\DeclareMathOperator{\frob}{\mathrm{Fr}}
\DeclareMathOperator{\Gal}{Gal}

\DeclareMathOperator{\Tr}{Tr}

\DeclareMathOperator{\swan}{Swan}
\DeclareMathOperator{\bb}{B}

\newcommand{\demi}{{\textstyle{\frac{1}{2}}}}

\newcommand{\sheaf}[1]{\mathcal{{#1}}}
\newcommand{\disj}{linearly disjoint}
\newcommand{\sheafm}[1]{\tilde{\sheaf{{#1}}}_{\ell}}

\DeclareMathSymbol{\gena}{\mathord}{letters}{"3C}
\DeclareMathSymbol{\genb}{\mathord}{letters}{"3E}

\def\sumb{\mathop{\sum \Bigl.^{\flat}}\limits}

\def\sums{\mathop{\sum \Bigl.^{*}}\limits}

\theoremstyle{plain}
\newtheorem{theorem}{Theorem}[section]

\newtheorem{lemma}[theorem]{Lemma}
\newtheorem{corollary}[theorem]{Corollary}

\newtheorem{proposition}[theorem]{Proposition}

\theoremstyle{remark}
\newtheorem{remark}[theorem]{Remark}

\theoremstyle{definition}

\newtheorem*{definition}{Definition}

\begin{document}

\title{The large sieve, monodromy and zeta functions of curves}
\author{E. Kowalski}
\address{Universit\'e Bordeaux I - A2X\\
351, cours de la Lib\'eration\\
33405 Talence Cedex\\
France}
\email{emmanuel.kowalski@math.u-bordeaux1.fr}
\begin{abstract}
We prove a large sieve statement for the average distribution of
Frobenius conjugacy classes in arithmetic monodromy groups over finite
fields. As a first application we prove a stronger version of a
result of Chavdarov on the ``generic'' irreducibility of the numerator
of the zeta functions in a family of curves with large monodromy. 
\end{abstract}

\maketitle

\section{Introduction}

In~\cite{chavdarov}, N. Chavdarov proves that, in an algebraic family
$C\ra U$ of smooth projective curves of genus $g$ over a finite field
$\Fp_q$, if the monodromy groups mod $\ell$ of the family are ``as large
as possible'' for almost all $\ell$, then the numerators
$\det(1-T\frob\,\mid\, H^1(\barre{C}_u,\Qq_{\ell}))$ of the zeta functions of
the curves $C_u$ of the family are ``almost all'' irreducible, and in
fact have splitting field as large as compatible with the existence of
the symplectic intersection pairing.
\par
Chavdarov's method, as sketched in the introduction
to~\cite{chavdarov}, is analogue in principle to the method used by
van der Waerden to show that ``most'' polynomials of given degree $d$
with integer coefficients have splitting field as large as
possible. This latter result was reproved in a simpler way and
stronger form by Gallagher~\cite{gallagher} using the large sieve
inequalities as a new analytic tool. This suggests trying to apply
similar ideas to Chavdarov's problem. In this paper, we show that this
is indeed possible, to some extent. This proof also yields a 
much stronger result than~\cite{chavdarov} in many cases; see
Theorem~\ref{th-chavdarov} and Theorem~\ref{th-chavdarov-2} for the
exact statements. The proof uses
some of the same tools, together with deep ideas of analytic number
theory and some new uniform estimates for $\ell$-adic Betti numbers
which may be of independent interest.
\par
\medskip
The plan of this paper is the following: in the next two sections, we
introduce the data involved and then state
our main bilinear form estimate 
from which we derive a  ``large sieve'' statement, essentially
in the same way as the classical case.  In the next two sections we
prove the bilinear form estimate. First,
Section~\ref{sec-betti} establishes some useful estimates for sums of
$\ell$-adic Betti numbers of ``Artin-like'' sheaves in various
circumstances (restricted unless the base space is a curve).
The proofs rely on the methods used by Katz in~\cite{katz-gsksm}
and~\cite{katz-betti} -- the difference being the parameters for
which uniformity is required. Then in Section~\ref{sec-bilinear},
we conclude the proof. 
\par
In the final sections we apply the sieve statement to prove our
form of Chavdarov's Theorem. The statements at least are accessible
(and of some interest) without knowledge of the techniques of étale
cohomology required for the proof. Families of a fixed genus
can be treated pretty quickly, but we expand some effort to obtain in
some cases a uniform result that can give information about curves of
genus $g$ over $\Fp_q$ when $q$ and $g$ are simultaneously large
(although $g$ must be much smaller than $q$). We also draw a few easy
consequences (Proposition~\ref{pr-sample1} and   
Proposition~\ref{pr-sample2}), as
illustrations of results which seem out of reach of Chavdarov's
method, but are is not meant as 
really important results in themselves.\footnote{\ For analytic number
  theorists, let us mention that, incidentally, we also get a
  version of Gallagher's 
estimate~\cite{gallagher} uniform in terms of the degree, see the
final remark in Section~\ref{sec-prelim-chav}.}
\par
One can hope that other applications of this method will arise, by
analogy with the situation 
in analytic number theory, where the ideas surrounding the large sieve
have been extremely successful since the original discovery by
Linnik. 
\par
\medskip
The related paper~\cite{weil-numbers} applies Chavdarov's Theorem and
some extra ingredients to the study of the characterization of abelian
varieties over finite fields or number fields by their torsion
fields. Also in~\cite{quad-twists}, we use the Betti number estimates
and a uniform Chebotarev density theorem to study the density of
quadratic twists of elliptic curves over function fields over finite
fields with rank $\geq 2$.
\par
\medskip
\textbf{Acknowledgments.} N. Katz found a serious mistake in the
first version of Section~\ref{sec-betti}; thanks for pointing it out
and for other enlightening remarks. 
\par
\medskip
\textbf{Background references.} Since we are using two important
themes in number theory which may not be equally known to the reader,
we mention a few general references. For the
large sieve, the reader may
consult~\cite{bombieri},~\cite{montgomery-survey}
or~\cite[Ch. 7]{ant}. For the approach to exponential sums via $\ell$-adic
methods and their applications (with which the author, for one, is not
so well acquainted), we
suggest~\cite[Ch. 9]{katz-sarnak}, \cite[Ch. 2,3]{katz-gsksm}, 
\cite[Ch. 1--3]{katz-se}, \cite[Sommes trig.]{sga4half}, or~\cite[\S
11.11]{ant}.
\par
\medskip
\textbf{Notation.} As usual, $|X|$ denotes the cardinality of a set,
$\mathfrak{S}_g$ is the symmetric group on $g$ letters. By $f\ll
g$ for $x\in X$, or $f=O(g)$ for $x\in X$, where $X$ is an
arbitrary set on which $f$ is defined, we mean synonymously that there
exists a constant $C\geq 0$ such that $|f(x)|\leq Cg(x)$ for all $x\in
X$. The ``implied constant'' is any admissible value of $C$. It may
depend on the set $X$ which is always specified or clear in context.

\section{Preliminaries}

Our main tool is a general estimate for a bilinear form
made up from representations of a system of lisse
$\barre{\Fp}_{\ell}$-sheaves on a variety over a finite field.
\par
The first basic data is therefore a base variety $U/\Fp_q$, where as
usual $\Fp_q$ denotes a finite field of characteristic $p$ with
$q$ elements. We assume that $U$ is smooth, affine, and geometrically
connected of dimension $d\geq 1$. 
\par
We denote by $\barre{\eta}$ the geometric generic point of $U$ and by
$\barre{U}$ the variety $U$ extended to $\barre{\Fp}_q$. We therefore
have the arithmetic fundamental group $\pi_1(U,\barre{\eta})$ and the
geometric fundamental group $\pi_1(\barre{U},\barre{\eta})$. Those fit
in an exact sequence
\begin{equation}\label{eq-piuns}
\begin{CD}
1 @>>> \pi_1(\barre{U},\barre{\eta}) @>>> \pi_1(U,\barre{\eta})
@>d>> \Gal(\barre{\Fp}_q/\Fp_q)\simeq \hat{\Zz} @>>>  1.
\end{CD}
\end{equation}
\par
For $n\geq 1$ and $u\in U(\Fp_{q^n})$, we denote by $\frob_{u,q^n}$
the geometric Frobenius automorphism at $u$ in
$\pi_1(U,\barre{\eta})$, i.e., the image of the inverse of the
canonical generator $x\mapsto x^{q^n}$ of the Galois group of
$\Fp_{q^n}$ via the map
$$
\Gal(\barre{\Fp}_{q^n}/\Fp_{q^n})\ra \pi_1(U,\barre{\eta})
$$
induced from the inclusion $\spec \Fp_{q^n}\ra
U$ which ``is'' $u$. In the above exact sequence we have then
$$
d(\frob_{u,q^n})=-n.
$$
\par
In most of our results, the base field (i.e., $q$) will be considered
fixed, althought the results will be uniform in $q$ so they can be
applied to $U\times \Fp_{q^n}$ for any $n\geq 1$. So most of the time
we just write $\frob_u$ instead of $\frob_{u,q}$ for $u\in U(\Fp_q)$.
\par
We also denote generically by $\frob$ the global
geometric Frobenius automorphism, acting for instance on $\ell$-adic
cohomology groups.
\par
\medskip
We now come to the sheaves on $U$ that we consider. We assume given a
set $\Lambda$ of primes $\not=p$, and for $\ell\in\Lambda$, a lisse
sheaf $\sheafm{F}$ of $\Fp_{\lambda}$-vector spaces of (fixed)
rank $r\geq 1$, where $\Fp_{\lambda}$ is a finite field of
characteristic $\ell$ (the degree of which over $\Fp_{\ell}$ may
depend on $\ell$). The basic example is when we have lisse sheaves
$\sheaf{F}_{\ell}$ of $\Zz_{\lambda}$-modules and 
$$
\sheafm{F}=\sheaf{F}_{\ell}/\ideal{m}_{\lambda}\sheaf{F}_{\ell},
$$
where $\Zz_{\lambda}$ is the ring of integers in a finite extension of
$\Qq_{\ell}$ with residue field $\Fp_{\lambda}$ and maximal ideal
$\ideal{m}_{\lambda}$. However, we do not assume that $\sheafm{F}$ is
of this type (of course, it will be in most applications).
\par
Equivalently (and this is the most convenient
viewpoint in terms of a first understanding at least),
$\sheafm{F}$ ``is'' a representation 
$$
\rho_{\ell}\,:\, \pi_1(U,\barre{\eta})\ra GL(r,\Fp_{\lambda}).
$$
From this description we can easily define the monodromy groups of
$\sheafm{F}$, or of $\rho_{\ell}$: the arithmetic monodromy
group $G_{\ell}\subset GL(r,\Fp_{\lambda})$ is the image of
$\rho_{\ell}$, and the geometric monodromy group $G^g_{\ell}$ is the
image of the subgroup $\pi_1(\barre{U},\barre{\eta})$. Thus
from~(\ref{eq-piuns}) we derive a commutative diagram
with exact rows and surjective downward arrows:
\begin{equation}\label{eq-cd}
\begin{CD}
1 @>>> \pi_1(\barre{U},\barre{\eta}) @>>> \pi_1(U,\barre{\eta})
@>d>> \hat{\Zz} @>>>  1\\
@. @VVV   @VVV  @V\varphi VV \\
1 @>>> G^g_{\ell} @>>>  G_{\ell}
@>m>> \Gamma_{\ell} @>>>  1,
\end{CD}
\end{equation}
where $\Gamma_{\ell}$ is a finite commutative (cyclic) group.
\par
\medskip
In the case where the sheaves $\sheafm{F}$ arise by reduction
of $\Zz_{\lambda}$-sheaves $\sheaf{F}_{\ell}$, as described
previously, one says that they form a \emph{compatible system} if for
every extension $\Fp_{q^n}/\Fp_q$, every $u\in U(\Fp_{q^n})$ and every
$\ell\in \Lambda$, the reversed characteristic polynomial of
$\frob_{u,q^n}$ acting on $\sheaf{F}_{\ell}$, i.e., the polynomial
$$
\det(1-T\frob_{u,q^n}\,\mid\, \sheaf{F}_{\ell})
$$
has coefficients in $\barre{\Qq}$ and is independent of $\ell$.
\par
\medskip
For any $\ell$ we will consider various sums involving irreducible
(complex valued) linear representations of $G_{\ell}$. For reasons
that will become clear during the proof of the main bilinear form
estimate (a certain phenomenon of ``imprimitivity''), we can not 
use all representations, but must ensure that those used are suitably
orthogonal when restricted to the subgroup $G^g_{\ell}$. 
\par
For this we have the following lemma.

\begin{lemma}\label{lm-ortho}
\emph{(1)} Let $\pi$, $\pi'$ be irreducible linear representations of
$G_{\ell}$. Then $\pi$ and $\pi'$ are equivalent when restricted to
$G_{\ell}^g$ if and only if there exists a character
$\psi\in\hat{\Gamma}_{\ell}$, the character group of
$\Gamma_{\ell}$, such that
$$
\pi=\pi'\otimes(\psi\circ m).
$$
\par
\emph{(2)} If $\pi$ and $\pi'$ are not equivalent restricted to
$G_{\ell}^g$, the representation
$\pi\otimes\tilde{\pi}'$ restricted to $G_{\ell}^g$ does not contain
the trivial representation, where $\tilde{\pi}'$ is the contragredient 
of $\pi'$. Otherwise it contains the trivial representation with
multiplicity equal to $|\hat{\Gamma}_{\ell}^{\pi}|$ where
$$
\hat{\Gamma}_{\ell}^{\pi}=\{\psi\in\hat{\Gamma}_{\ell}\,\mid\,
\pi\simeq \pi\otimes(\psi\circ m)\}.
$$
\end{lemma}

\begin{proof}
We will identify characters $\psi\in \hat{\Gamma}_{\ell}$ with
characters of $G_{\ell}$ by $\psi(x)=\psi(m(x))$ for $x\in G_{\ell}$.
\par
For any representation $\tau$ of $G_{\ell}$, let $j(\tau)$ denote the
multiplicity of the trivial representation in the 
restriction of $\tau$ to $G_{\ell}^g$. This is given by
\begin{align*}
j(\tau)&=\frac{1}{|G_{\ell}^g|}\sum_{x\in G_{\ell}^g}{\Tr\tau(x)}\\
&=\frac{1}{|G_{\ell}^g|}
\sum_{x\in G_{\ell}}{\frac{1}{|\Gamma_{\ell}|}
\sum_{\psi\in\hat{\Gamma}_{\ell}}{\psi(m(x))}\Tr\tau(x)}\\
&=
\sum_{\psi\in\hat{\Gamma}_{\ell}}{\frac{1}{|G_{\ell}|}
\sum_{x\in G_{\ell}}{\psi(x)\Tr \tau(x)}}
\end{align*}
(by orthogonality of characters of $\Gamma_{\ell}$),
which is the sum of the multiplicities of the characters $\psi$ of
$G_{\ell}$ in $\tau$. (This interpretation being of course also
available by a simple application of Frobenius reciprocity).
\par
Applying this to $\tau=\pi\otimes\tilde{\pi}'$, it follows that $1$ is
contained in the restriction of $\tau$ to $G_{\ell}^g$ if and only if
there exists a $\psi$ such that $\psi$ is contained in
$\pi\otimes\tilde{\pi}'$. However if that is the case, the trivial
representation is contained in $\pi\otimes \tilde{\pi}'\barre{\psi}$,
but as $\pi$ and $\pi'\psi$ are irreducible, this
means that $\pi\simeq \pi'\otimes\psi$ as representations of
$G_{\ell}$. This shows the ``only if'' part of the lemma, and the
other direction is trivial since $\psi$ restricts to the trivial
character of $G_{\ell}^g$.
\par
The first part of (2) is contained in the previous paragraph. The
assertion about the multiplicity is also clear: if
$\pi'=\pi\otimes\psi_0$ is a twist of $\pi$ by a character $\psi_0$, the
multiplicity $j(\pi\otimes\tilde{\pi}')$ is the sum of multiplicities
of the characters $\psi$ in $\pi\otimes\tilde{\pi}\barre{\psi}_0$,
each of which is equal to $1$ if $\pi\simeq \pi\otimes(\psi\psi_0)$,
and $0$ otherwise, i.e, it is equal to $1$ if $\psi\psi_0\in
\hat{\Gamma}_{\ell}^{\pi}$ and $0$ otherwise. So the total
multiplicity is the number of elements in
$\hat{\Gamma}_{\ell}^{\pi}$. 
\end{proof}

\begin{remark}
If $\pi\simeq \pi\otimes \psi$, we must have $\psi(x)=1$ whenever
$\Tr\pi(x)\not=0$. But if $\deg(\pi)>1$, it is well-known that there
are elements $x\in G_{\ell}$ with $\Tr\pi(x)=0$, and then the value of
$\psi$ is not determined. Of course, ``in general'', we have
$\hat{\Gamma}_{\ell}^{\pi}=1$, but (for instance), for any
representation $\pi$ of degree $2$ of a dihedral group $D_n$, $n$ even (of
order $2n$), there is a character $\psi$ with $\pi\simeq
\pi\otimes\psi$. 
\end{remark}

Say that two representations $G_{\ell}$ are geometrically equivalent
if their restrictions to $G_{\ell}^g$ are equivalent, or (by the
lemma) if and only if they differ by a twist by a character of
$\Gamma_{\ell}$. 
We now assume chosen a set  $\Pi_{\ell}$ of representatives 
of the irreducible representations of $G_{\ell}$ for this equivalence
relation. Using these and characters of $\Gamma_{\ell}$ one can
parameterize all irreducible representations of $G_{\ell}$ as follows: 
they are of the form $\pi\otimes\psi$ where $\pi\in
\Pi_{\ell}$ and $\psi\in\hat{\Gamma}_{\ell}$; the representation $\pi$
is unique, but $\psi$ is only unique up to multiplication by an
element of the group $\hat{\Gamma}_{\ell}^{\pi}$ defined in (2) of the
previous lemma.
\par
This ambiguity requires us to control the size of those groups
$\hat{\Gamma}_{\ell}^{\pi}$. We will assume that for all
$\ell\in\Lambda$ and $\pi\in \Pi_{\ell}$, we have
\begin{equation}\label{eq-kappa}
|\hat{\Gamma}_{\ell}^{\pi}|\leq \kappa
\end{equation}
for some fixed $\kappa\geq 1$.
\par
Here are useful cases when we can get such a bound.
\begin{lemma}\label{lm-kappa}
\emph{(1)} Assume that for all $\ell$ we have
$G^g_{\ell}=SL(r,\Fp_{\ell})$.
Then~\emph{(\ref{eq-kappa})} holds with $\kappa=r$.
\par
\emph{(2)} Assume that $r$ is even and that for all $\ell$ we have
$G^g_{\ell}=Sp(r,\Fp_{\ell})$, the symplectic group for some
non-degenerate alternating form $\langle \cdot,\cdot\rangle$ on
$\Fp_{\ell}^r$, and that $G_{\ell}$ 
is a subgroup of the group $SSp(r,\Fp_{\ell})$ of symplectic
similitudes, i.e., for $g\in G_{\ell}$ we have $\langle
gv,gw\rangle=m(g)\langle v,w\rangle$ for some
$m(g)\in\Fp_{\ell}^{\times}$, called the multiplicator of $g$. 
Then~\emph{(\ref{eq-kappa})} holds with $\kappa=2$.
\end{lemma}

\begin{proof}
(1) If $\pi$ is an irreducible representation of $G_{\ell}$ and
$\psi\in\hat{\Gamma}_{\ell}^{\pi}$, then $\psi$ is trivial on the
center $Z_{\ell}$ of $G_{\ell}$. For any $x\in G_{\ell}$, we can write
$$
x^r=(\det x) y
$$
with $y\in SL(r,\Fp_{\ell})=G_{\ell}^g$. Hence $\det(x)\in
\Fp_{\ell}^{\times}\cap G_{\ell}\subset Z_{\ell}$ and therefore
$\psi(\det(x))=1$, and $\psi(x^r)=\psi(\det(x))\psi(y)=1$. So $\psi$
is of order at most $r$, and since it is a character of a cyclic
group, there are at most $r$ such characters, giving~(\ref{eq-kappa})
with $\kappa=r$.
\par
(2) The argument is similar except that we now have $x^2=m(x)y$ with
$y\in Sp(r,\Fp_{\ell})$ (since $m(ax)=a^2 x$ for scalar
$a$), so $\psi(x)^2=1$.
\end{proof}

We will often simply write (when $\ell$ is cleary specified, e.g. as a
summation parameter occuring before)
$$
\sums_{\pi}{\alpha(\pi,\ell,\ldots)},\quad 
\sums_{\pi\not=1}{\alpha(\pi,\ell,\ldots)}
$$
for, respectively, a sum over all the irreducible representations $\pi\in
\Pi_{\ell}$ of $G^g_{\ell}$ or for a sum over all those which are
non-trivial on $G_{\ell}^g$. Similarly, a sum of the type
$$
\sums_{\pi}{\sum_{\psi}{\alpha(\pi,\psi\ldots)}}
$$
means (unless otherwise specified) that $\pi\in\Pi_{\ell}$ and
$\psi\in\hat{\Gamma}_{\ell}/\hat{\Gamma}_{\ell}^{\pi}$; in other
words, this is a sum over all irreducible representations of
$G_{\ell}$, parameterized as described previously.
\par
\medskip
We need various estimates involving sums of dimensions of the
representations in $\Pi_{\ell}$. We will phrase them in terms of
upper bounds for the ``dimensions'' of $G_{\ell}$ and of the set 
$G_{\ell}^{\sharp}$ of its conjugacy classes: let $s$ and $t$ be such
that the inequalities
\begin{equation}\label{eq-s}
|G_{\ell}|\leq c_1\ell^s,\quad |G^{\sharp}_{\ell}|\leq c_2\ell^t
\end{equation}
hold for all primes $\ell\in\Lambda$, $c_1$ and $c_2$ being two
given constants. Note that of course $s=t=r^2$ is always possible with
$c_1=c_2=1$ (and that in fact this does not in general significantly
affect the applications).

\begin{lemma}\label{lm-sizes}
\emph{(1)} We have
$$
\sum_{\pi\in\Pi_{\ell}}{\dim \pi}\leq (c_2\ell^{s+t})^{1/2},
$$
and for all $\pi\in\Pi_{\ell}$ we have
$$
\dim \pi\leq (c_1\ell^s)^{1/2}.
$$
\par
\emph{(2)} If $G_{\ell}^g=SL(r,\Fp_{\ell})$, the
estimates~\emph{(\ref{eq-s})} hold with
$$
c_1=1,\quad s=r^2,\quad c_2=6^{r},\quad t=r.
$$
\par
\emph{(3)} If $r$ is even, $G_{\ell}^g=Sp(r,\Fp_{\ell})$ and
$G_{\ell}\subset SSp(r,\Fp_{\ell})$,  the
estimates~\emph{(\ref{eq-s})} hold with
$$
c_1=1,\quad s=1+\frac{r(r+1)}{2},\quad c_2=6^{r/2},\quad t=r/2+1.
$$
\end{lemma}

\begin{proof}
For a representation of a finite group $G$, the dimension is always
$\leq |G|^{1/2}$, and the sum of the dimension is bounded by Cauchy's
inequality by
$$
\sum_{\pi}{\dim \pi}\leq |G^{\sharp}|^{1/2}|G|^{1/2},
$$
so that (1) is a direct translation of~(\ref{eq-s}).
\par
(2) is obvious, noticing that the number of conjugacy classes in
$G_{\ell}$ is at most 
$$
|\Gamma_{\ell}||G_{\ell}^{g,\sharp}|\leq
(\ell-1)||G_{\ell}^{g,\sharp}|
\leq (\ell-1)(6\ell)^{r-1}
$$
(by~\cite[Lemma 1.4]{liebeck-pyber} for instance; the factor $6^{g-1}$
takes into account the non-semisimple conjugacy classes).
\par
(3) is similar using the formula for the cardinality of
$Sp(2g,\Fp_{\ell})$, and~\cite[Lemmas
1.3,1.6]{liebeck-pyber} for the conjugacy classes.
\end{proof}

\par
\medskip
Our last definition is also of crucial importance for the bilinear
form estimate:  
\begin{definition}\label{def-disj-twist}
We say that the family
$(\sheafm{F})$ is \emph{\disj} if for all $\ell$ and $\ell'$ in
$\Lambda$, with $\ell\not=\ell'$, the product map
$$
\pi_1(\barre{U},\barre{\eta})\ra G^g_{\ell}\times G^g_{\ell'}
$$
is surjective.
\end{definition}

This is a fairly natural
independence notion for the various monodromy groups. In many cases it
will hold for group-theoretical reasons simply because the $G^g_{\ell}$ are
``large'' groups and close to simple; this is related to Goursat's
lemma, and we quote here the version in~\cite[Pr. 5.1]{chavdarov},
(specialized for $2$ factors):

\begin{lemma}\label{lm-goursat}
Let $G_1$ and $G_2$ be finite groups such that every normal subgroup
of $G_i$ is contained in the center $C_i$, and such that $G_1/C_1$ and
$G_2/C_2$ are distinct, simple and non-abelian. Then no proper subgroup
$G\subset G_1\times G_2$ projects surjectively on both $G_1$ and $G_2$.
\end{lemma}

This is typically applied with $G_1=G^g_{\ell}$, $G_2=G^g_{\ell'}$
and $G$ the image of $\pi_1(\barre{U},\barre{\eta})\ra G_1\times G_2$
which \emph{does} project surjectively on both factors.
\par
For instance, this shows:
\begin{corollary}\label{cor-goursat}
\emph{(1)}  Let $r$ be even and let $(\sheafm{F})$ be a family of
sheaves as above such that $G^g_{\ell}=Sp(r,\Fp_{\ell})$ for all
$\ell$ in $\Lambda$, with $\ell\geq 5$ if $r=2$ and $\ell\geq 3$ if
$r=4$. Then the family is \disj .
\par
\emph{(2)} Let $(\sheafm{F})$ be a family of sheaves as above such that
$G^g_{\ell}=SL(r,\Fp_{\ell})$ for all $\ell$ in $\Lambda$, with
$\ell\geq 5$ if $r=2$. Then the family is \disj . 
\end{corollary}

This follows because it is very classical that the center of
$Sp(r,\Fp_{\ell})$ (resp. $SL(r,\Fp_{\ell})$) is $\pm 1$ (and is the
only non-trivial normal subgroup) and $Sp(r,\Fp_{\ell})/\{\pm 1\}$
(resp. $SL(r,\Fp_{\ell})/\{\pm 1\}$) is a simple non-abelian 
group in the cases described (see e.g.~\cite[Th. 5.1,
Th. 4.9]{artin}). On the other hand,  
notice the lemma can not be applied for orthogonal groups. (For
instance, if $\ell$, $\ell'$ are odd, the proper subgroup
$$
\{(x,y)\in O(r,\Fp_{\ell})\times O(r,\Fp_{\ell'})\,\mid\,
\det(x)=\det(y)\},
$$
where the equality makes sense because the determinants are $\pm 1$,
clearly projects surjectively onto both factors).

\section{Bilinear form estimates and large sieve for algebraic families}
\label{sec-main}

We now state the bilinear form estimate which is our main tool.

\begin{theorem}\label{th-main}
Let $U$ be a variety and $(\sheafm{F})$ a family of
sheaves as above, with given sets $\Pi_{\ell}$ of irreducible
representations which are representatives for geometric
equivalence. Assume that the family is \disj, that it
satisfies~\emph{(\ref{eq-kappa})} 
and moreover that $U$ and $(\sheafm{F})$ satisfy one of the
following conditions:
\par
\quad \emph{(i)} $U$ is a smooth affine curve and $(\sheafm{F})$ arises
from a compatible system of integral $\ell$-adic sheaves;
\par
\quad \emph{(ii)} For all $\ell\in\Lambda$, the order of $G^g_{\ell}$ is
prime to $p$.
\par
Then there exists constants $C\geq 0$ and $A\geq 0$ such that we have
\begin{equation}\label{eq-main-bil}
\sum_{\ell\leq L}{\sums_{\pi\not=1}{
\Bigl|
\sum_{u\in U(\Fp_q)}{\alpha(u)\Tr(\pi\circ\rho_{\ell})(\frob_u)}
\Bigr|^2
}}
\leq (\kappa q^d+Cq^{d-1/2}L^{A})
\sum_{u\in U(\Fp_q)}{|\alpha(u)|^2},
\end{equation}
for any $L\geq 1$ and any complex coefficients
$\alpha(u)$.  
\par
In case \emph{(i)}, we can take $A=1+s+t/2$, and the 
constant $C$ depends only on $\barre{U}$, the ``geometric'' compatible
system $(\sheaf{F}_{\ell})$ on $\barre{U}$ and the constants $c_1$ and
$c_2$.  In case~\emph{(ii)} we can take $A=1+5s/2+t/2$, and the 
constant $C$ depends only on $\barre{U}$, $c_1$ and $c_2$.
\par
In particular the estimate can be applied uniformly for
$U\otimes\Fp_{q^n}$ for any $n\geq 1$.  
\end{theorem}

Note that the left-hand side of~(\ref{eq-main-bil}) is in fact
independent of the choice of representative sets $\Pi_{\ell}$.

\begin{remark}
Here are a few remarks, most of which are of a general nature and 
are standard observations for any type of bilinear form estimate. 
\par
(1) The estimate~(\ref{eq-main-bil}) is most interesting when $L$ is
small enough that $L^A\leq q^{1/2}$, so that the sum of the two terms
$q^d$ and $q^{d-1/2}L^A$ is still of size $q^d$, which is roughly the
number of terms in the inner sum over $u\in U(\Fp_q)$ by the Lang-Weil
estimate $|U(\Fp_q)|=q^d+O(q^{d-1/2})$. 
\par
(2) The restriction to $\pi\not=1$ in the summation
in~(\ref{eq-main-bil}) is essential: the additional contribution of
the trivial representations would give the quadratic form
$$
\sum_{\ell\leq L}{\Bigl|\sum_{u\in U(\Fp_q)}{\alpha(u)}\Bigr|^2}
$$
which by Cauchy's inequality has norm $|U(\Fp_q)|L\asymp q^dL$, which
exceeds $(\kappa q^d+q^{d-1/2}L^A)$ in the most interesting ranges
where $L$ is small as in the previous remark.
\par
(3) In~(\ref{eq-main-bil}), the trivial bound has $(\kappa
q^d+q^{d-1/2}L^A)$ replaced by
$$
|U(\Fp_q)|\sum_{\ell\leq L}{\sum_{\pi\not=1}{1}}\asymp
q^dL^{1+t}
$$
(i.e., the ratio is bounded from above and below; we assume that $t$
is chosen optimally).
On the other hand, from general principles (see e.g.~\cite[\S 7]{ant}), the
best possible result is essentially with
$$
|U(\Fp_q)|+\sum_{\ell\leq L}{\sum_{\pi\not=1}{1}}\asymp q^d+L^{1+t},
$$
and nowadays it is usually estimates of similar strength which are
called large sieve inequalities, even when no connection with sieve
theory exists. 
\par
Thus from the point of view of the study of bilinear forms in
modern analytic number theory, the estimate~(\ref{eq-main-bil}) is
much too weak to deserve the name of large sieve. However, we
\emph{are} using it mainly to derive a sieve-type result which
corresponds to Linnik's original description of a ``large'' sieve, so
we use the word in this sense. 
\par
(4) Using geometric considerations one can get similar results for
more general $U$, for instance by dealing with irreducible components
one by one. Or if $U$ is the disjoint
union of $U_1$ and $U_2$, with $U_1$ a dense open subscheme which is
smooth affine and geometrically connected, and $U_2$ closed of
codimension $\geq 1$, the sum on the left of~(\ref{eq-main-bil}) is at
most twice the sum
$$
\Bigl|
\sum_{u\in U_1(\Fp_q)}{\alpha(u)\Tr(\pi\circ\rho_{\ell})(\frob_u)}
\Bigr|^2
+\Bigl|
\sum_{u\in U_2(\Fp_q)}{\alpha(u)\Tr(\pi\circ\rho_{\ell})(\frob_u)}
\Bigr|^2
$$
and~(\ref{eq-main-bil}) applies to the first sum while the second
has a contribution $\ll mq^{d-1}L^{1+t}$ by Remark~(3), where $m$ is
the number of irreducible components of
$U_2\otimes\barre{\Fp}_q$. Another case would be to have a map $U\ra V$
with ``most'' fibers being smooth, affine, connected curves on which
the induced sheaves have the same monodromy as on $U$, and for
which the constant $C$ in~(\ref{eq-main-bil}) happens to be uniformly
bounded for all such fibers.
\par
(5) Formula~(\ref{eq-delta-final}) below gives a more explicit bound
which may be better in some cases where more is known about the
$G_{\ell}$ (although it's not clear how much of a difference it
would make in applications), for instance the maximal dimension of an
irreducible representation.\footnote{\ This is certainly known for the
groups $Sp(2g,\Fp_{\ell})$ that we will use below, but the author
doesn't know where to find it...}
\par
(6) Finally, we note that a standard heuristic understanding of the
strength of the classical large sieve inequality (see
e.g.~\cite[7.5]{ant} for a proof)
$$
\sum_{q\leq Q}{\sums_{\chi\mods{q}}{\Bigl|
\sum_{n\leq N}{a_n\chi(n)}
\Bigr|^2}}\leq (N-1+Q^2)\sum_{n\leq N}{|a_n|^2}
$$
(where $\sums$ indicates a sum over \emph{primitive} Dirichlet
characters modulo $q$), is that in its range of effectiveness (i.e.,
$Q^2\leq N$) it is as strong as the Generalized Riemann Hypothesis, as
it gives for instance
$$
\sum_{n\leq N}{\mu(n)\chi(n)}\ll \sqrt{N}
$$
on average. And indeed this inequality is used as a substitute for GRH
in many applications. In 
view of this and the fact that we already know the Riemann Hypothesis
over finite fields by Deligne's work, one may think that a large sieve
inequality would be either trivial to prove or without application (or
both) in this context. That it is not the case illustrates two points:
first, that in the case of sums over a variety of dimension $\geq 2$,
and for $L^A\ll \sqrt{q}$, the inequality~(\ref{eq-main-bil}) gives on
average \emph{square-root cancellation} in
$$
\sum_{u\in U(\Fp_q)}{\alpha(u)\Tr(\pi\circ\rho_{\ell})(\frob_u)},
$$
which is \emph{stronger} than the Riemann Hypothesis because the
latter only provides a saving of $q^{-1/2}$ in general, not
$q^{-d/2}$. Thus~(\ref{eq-main-bil}) contains information about the
finer average distribution of the zeros of the $L$-functions that are involved
in those exponential sums. The second point, which is valid even for
$d=1$, is that the large sieve inequality is not only about
cancellation, but about \emph{uniformity} in estimates. Hence it is
not surprising that the ``only'' difficulty in
proving~(\ref{eq-main-bil}), from the Riemann Hypothesis, is a
question of uniform bounds for the error terms coming after applying
Deligne's results.
\end{remark}

We will prove Theorem~\ref{th-main} in
Section~\ref{sec-bilinear}. For the moment, we derive a large sieve
estimate concerning the average distribution of the Frobenius
conjugacy classes in $G_{\ell}$. 
\par
Let $L\geq 2$ and 
suppose that for $\ell\in\Lambda$, $\ell\leq L$, we select some
conjugacy-invariant subset $\Omega(\ell)$ of 
$G_{\ell}$ with cardinality $\omega(\ell)$, such that 
$$
m(x)=\varphi(-1)\in\Gamma
$$
for all $x$ and $\ell$ (where $m\,:\, G_{\ell}\ra \Gamma$ and
$\varphi$ are defined by the
commutative diagram~(\ref{eq-cd}); recall that
$d(\frob_u)=-1\in\hat{\Zz}$ for $u\in U(\Fp_q)$).
\par
Let then
$$
P(u,L)=\sum_{\stacksum{\ell\leq L}{\rho_{\ell}(\frob_u)\in \Omega(\ell)}}{1}
$$
for $u\in U(\Fp_q)$ and
$$
P(L)=\sum_{\ell\leq L}{\omega(\ell)|G^g_{\ell}|^{-1}}.
$$
The large sieve statement says that for ``most'' $u$, the value of
$P(u,L)$ is close to the average value $P(L)$, this being measured by
the variance.

\begin{proposition}\label{pr-large-sieve}
With $U$ and $(\sheafm{F})$ satisfying one of the
assumptions of Theorem~\ref{th-main}, 
we have 
\begin{equation}\label{eq-large-sieve}
\sum_{u\in U(\Fp_q)}{(P(u,L)-P(L))^2}\leq
(\kappa q^d+Cq^{d-1/2}L^{A})P(L),
\end{equation}
where the constants $C$ and $A$ are the same as in
Theorem~\ref{th-main}. In particular, the cardinality of the sifted
set 
$$
S(U,\Omega;L)=\{
u\in U(\Fp_{q})\,\mid\,
\frob_u\notin\Omega(\ell)\text{ for all } \ell\leq L
\}
$$
satisfies
\begin{equation}\label{eq-sieve}
|S(U,\Omega;L)|\leq (\kappa q^d+Cq^{d-1/2}L^{A})P(L)^{-1}.
\end{equation}
\end{proposition}

\begin{proof}
First~(\ref{eq-sieve}) follows trivially from~(\ref{eq-large-sieve})
since $P(u,L)=0$ for $u\in S(U,\Omega;L)$, so that the left-hand side
of the latter inequality is at least equal to $P(L)^2|S(U,\Omega;L)|$.
\par
So we prove~(\ref{eq-large-sieve}); the argument is in large part a
jazzed-up version of the one in~\cite{gallagher}.
Let $\chi_{\ell}$ be
the characteristic function of $\Omega(\ell)$. Since $\Omega(\ell)$ is
invariant by conjugation, we can expand it in Fourier
series using the representations of $G_{\ell}$. Using the
parameterization as $\pi\otimes\psi$ with $\pi\in
\Pi_{\ell}$ and $\psi\in
\hat{\Gamma}_{\ell}/\hat{\Gamma}_{\ell}^{\pi}$, we can write this
expansion as
\begin{equation}\label{eq-fourier-exp}
\chi_{\ell}(x)=\sum_{\pi\in\Pi_{\ell}}{\sum_{\psi\in\hat{\Gamma}_{\ell}/
\hat{\Gamma}_{\ell}^{\pi}}
{\hat{\chi}_{\ell}(\psi,\pi)
\psi(x)\Tr \pi(x)}}
\end{equation}
with 
\begin{equation}\label{eq-fourier-coeff}
\hat{\chi}_{\ell}(\psi,\pi)=\frac{1}{|G_{\ell}|}\sum_{x\in \Omega(\ell)}
{\barre{\psi(x)}\overline{\Tr\pi(x)}}=\frac{1}{|\Gamma_{\ell}|}
\psi(\varphi(1))\gamma(\pi),
\end{equation}
where
$$
\gamma(\pi)=\frac{1}{|G^g_{\ell}|}\sum_{x\in \Omega(\ell)}
{\overline{\Tr\pi(x)}}.
$$
Thus we have
\begin{equation}\label{eq-harm-princ}
\gamma(1)=\omega(\ell)|G^g_{\ell}|^{-1}.
\end{equation}
Also by orthonormality of the characters of $G_{\ell}$ we have
$$
\frac{\omega(\ell)}{|G_{\ell}|}
=\frac{1}{|G_{\ell}|}\sum_{x\in G_{\ell}}{
|\chi_{\ell}(x)|^2}
=
\sums_{\pi}{\sum_{\psi}{|\hat{\chi}_{\ell}(\psi,\pi)|^2}}
=\sums_{\pi}{\sum_{\psi}{\frac{1}{|\Gamma_{\ell}|^2}
|\gamma(\pi)|^2}}
=\frac{1}{|\Gamma_{\ell}|}\sums_{\pi}{\frac{|\gamma(\pi)|^2}
{|\hat{\Gamma}_{\ell}^{\pi}|}},
$$
hence
\begin{equation}\label{eq-positivity}
\sums_{\pi\not=1}{\frac{|\gamma(\pi)|^2}{|\hat{\Gamma}_{\ell}^{\pi}|^2}}
\leq 
\sums_{\pi\not=1}{\frac{|\gamma(\pi)|^2}{|\hat{\Gamma}_{\ell}^{\pi}|}}
\leq
\sums_{\pi}{\frac{|\gamma(\pi)|^2}{|\hat{\Gamma}_{\ell}^{\pi}|}}
=
\frac{\omega(\ell)}{|G^g_{\ell}|}.
\end{equation}
\par
By~(\ref{eq-fourier-exp}),~(\ref{eq-fourier-coeff}) and the fact that
$\psi(\rho_{\ell}(\frob_u))=\psi(\varphi(-1))$, we get that
for $u\in U(\Fp_q)$ and $\ell\leq L$ we have
\begin{align*}
\chi_{\ell}(\rho_{\ell}(\frob_u))&=
\sums_{\pi}{\sum_{\psi}
{\frac{1}{|\Gamma_{\ell}|}\psi(\varphi(1))\gamma(\pi)
\psi(\varphi(-1))\Tr (\pi\circ\rho_{\ell})(\frob_u)}}\\
&=\sums_{\pi}{\gamma(\pi)\Tr (\pi\circ\rho_{\ell})(\frob_u)
\Bigl(\frac{1}{|\Gamma_{\ell}|}\sum_{\psi}{\psi(\varphi(1))
\psi(\varphi(-1))}\Bigr)}\\
&=\sums_{\pi}{\frac{\gamma(\pi)}{|\hat{\Gamma}_{\ell}^{\pi}|}
\Tr (\pi\circ\rho_{\ell})(\frob_u)}
\end{align*}
hence (since $\hat{\Gamma}_{\ell}^1=1$)
\begin{equation}\label{eq-pay}
P(u,L)=\sum_{\ell\leq L}{\sums_{\pi}{
\frac{\gamma(\pi)}{|\hat{\Gamma}_{\ell}^{\pi}|}
\Tr(\pi\circ\rho_{\ell})(\frob_u)}}=
P(L)+\sum_{\ell\leq L}{\sums_{\pi\not=1}{
\frac{\gamma(\pi)}{|\hat{\Gamma}_{\ell}^{\pi}|}
\Tr(\pi\circ\rho_{\ell})(\frob_u)}}
\end{equation}
using~(\ref{eq-harm-princ}).
\par
Denote by $R(u,L)$ the second term on the right-hand side (the sum over
$\ell\leq L$ and $\pi\not=1$). By Cauchy's inequality
and~(\ref{eq-positivity}) we have  
\begin{align*}
\sum_{u\in U(\Fp_q)}{|R(u,L)|^2}&=
\sum_{\ell\leq L}{\sums_{\pi\not=1}{
\frac{\gamma(\pi)}{|\hat{\Gamma}_{\ell}^{\pi}|}
\sum_{u\in U(\Fp_q)}{R(u,L)\Tr(\pi\circ\rho_{\ell})(\frob_u)}}}\\
&\leq
\Bigl(
\sum_{\ell\leq L}{\sums_{\pi\not=1}{
\frac{|\gamma(\pi)|^2}{|\hat{\Gamma}_{\ell}^{\pi}|^2}}}
\Bigr)^{1/2}
\Bigl(
\sum_{\ell\leq L}{\sums_{\pi\not=1}{\Bigl|
\sum_{u\in U(\Fp_q)}{R(u,L)\Tr(\pi\circ\rho_{\ell})(\frob_u)}
\Bigr|^2}}
\Bigr)^{1/2}\\
&\leq P(L)^{1/2}
\Bigl(
\sum_{\ell\leq L}{\sums_{\pi\not=1}{\Bigl|
\sum_{u\in U(\Fp_q)}{R(u,L)\Tr(\pi\circ\rho_{\ell})(\frob_u)}
\Bigr|^2}}
\Bigr)^{1/2}.
\end{align*}
We can apply Theorem~\ref{th-main} to the last sum, getting (after squaring)
$$
\Bigl(\sum_{u\in U(\Fp_q)}{|R(u,L)|^2}\Bigr)^2\leq
P(L)(\kappa q^d+Cq^{d-1/2}L^{A})
\sum_{u\in U(\Fp_q)}{|R(u,L)|^2}
$$
so that
$$
\sum_{u\in U(\Fp_q)}{|R(u,L)|^2}\leq (\kappa q^d+Cq^{d-1/2}L^{A})P(L),
$$
which concludes the proof since by~(\ref{eq-pay}) we have
$$
\sum_{u\in U(\Fp_q)}{(P(u,L)-P(L))^2}=
\sum_{u\in U(\Fp_q)}{|R(u,L)|^2}.
$$
\end{proof}

\section{Estimates for sums of Betti numbers}
\label{sec-betti}

In this section we will prove some estimates for Betti numbers of
$\ell$-adic sheaves needed in the proof of the main estimate in
Section~\ref{sec-bilinear}. 
\par
For a separated scheme of finite type $U$ over $\barre{\Fp}_q$ of
dimension $d\geq 1$ and any prime
$\ell\not=p$ we denote as usual
\begin{gather*}
h^i_c(U,\sheaf{F})=\dim H^i_c(U,\sheaf{F}),\quad
h^i(U,\sheaf{F})=\dim H^i(U,\sheaf{F}),\\
\sigma_c(U,\sheaf{F})=\sum_{i}{h^i_c(U,\sheaf{F})},\quad
\sigma(U,\sheaf{F})=\sum_{i}{h^i(U,\sheaf{F})},\\
\chi_c(U,\sheaf{F})=\sum_{i}{(-1)^ih^i_c(U,\sheaf{F})},\quad
\chi(U,\sheaf{F})=\sum_{i}{(-1)^ih^i_c(U,\sheaf{F})},
\end{gather*}
where $\sheaf{F}$ can be either a $\barre{\Qq}_{\ell}$-sheaf on $U$
or an $\barre{\Fp}_{\ell}$-sheaf. We also write
$$
\sigma'_c(U,\sheaf{F})=\sum_{i<2d}{h^i_c(U,\sheaf{F})}\leq
\sigma_c(U,\sheaf{F})
$$
for the sum of all Betti numbers except the topmost one. 
\par
\medskip
We first consider the case of curves. Here the situation will be as
follows: $U/\barre{\Fp}_q$ is a smooth affine connected curve,
$\rho\,:\, \pi_1(U,\barre{\eta})\ra G$ is a surjective 
group homomorphism with $G$ finite, and $\pi$ is a representation of
$G$, with values in some $\barre{\Qq}_{\ell}$-vector space of finite
dimension, for some $\ell\not=p$.
We can form the composite $\pi\circ 
\rho$ to obtain a lisse $\ell$-adic sheaf on $U$, which is denoted
$\pi(\rho)$. Then 
we wish to find bounds for the sum of Betti numbers
$\sigma'_c(U,\pi(\rho))$ which are polynomial in the size of $G$
(or the degree $\dim\pi$ of $\pi$).
\par
We  do this for
$\rho$ of a special type, which we describe in a slightly more general
case than will be needed in the next section: $G$ is a product
$$
G=\prod_{1\leq i\leq k}{G_i}
$$
where $G_i$ is a subgroup of $GL(r,\Fp_{\lambda_i})$ for $1\leq i\leq
k$, $\lambda_i$ a power of a prime $\ell_i\not=p$ (the $\ell_i$ are
not necessarily distinct),  
and the representation $\rho$ is a tensor product
$\rho_1\otimes\cdots\otimes\rho_k$ where the $\rho_i$ correspond to
lisse sheaves $\sheaf{F}_i/\ell_i\sheaf{F}_i$, which are the
reductions modulo $\ell_i$ of sheaves
$\sheaf{F}_i$ of $\Zz_{\lambda_i}$-modules which are part of a
compatible system $(\sheaf{F}_{\ell})$. Here $\Zz_{\lambda_i}$ is
the ring of integers of a finite extension $\Qq_{\lambda_i}$ of
$\Qq_{\ell_i}$ with residue field $\Fp_{\lambda_i}$.
\par
Our goal is:

\begin{proposition}\label{pr-boundcurve}
With notation as above, we have
the bound
$$
\sigma'_c(U,\pi(\rho))\leq 
C(U,(\sheaf{F}_i),k)(\dim \pi),
$$
for some constant $C(U,(\sheaf{F}_i),k)$ depending only on $U$, $k$
and the compatible system, but not on $\pi$. One can take
\begin{equation}\label{eq-ww}
C=1-\chi_c(U,\Qq_{\ell})+kw
\end{equation}
where $w\geq 0$ is the sum of the Swan conductors of all $\sheaf{F}_i$ at
the points at infinity for $U$, as described below, which is
independent of $i$.
\end{proposition}

We start by recalling and setting up the description of the
ramification structure of sheaves on $U$, as described for instance
in~\cite[Ch. 1]{katz-gsksm}. Let $C$ be the smooth
projective model of $U$ and $S=C-U$ the non-empty finite set of
``points at infinity''. Let $M$ be a $\Zz[1/p]$-module on which
$\pi_1(U,\barre{\eta})$ acts through a finite discrete quotient. For
each point $x\in S$, there is a certain direct sum decomposition of
$M$ seen as representation of the inertia group $I_x$ at $x$ of the type
$$
M=\bigoplus_{t\geq 0}{M_x(t)}
$$
where each $M_x(t)$ is $I_x$-stable. All but
finitely many of the $M_x(t)$ vanish, and those $t$ for which
$M_x(t)\not=0$ are called the breaks of $M$ at $x$. If $M$ is free
over some $\Zz[1/p]$-algebra $A$ (e.g., $A=\Fp_{\ell}$, $\Zz_{\ell}$ or
$\Qq_{\ell}$), the Swan
conductor of $M$ at $x$ is then defined by 
$$
\swan_x(M)=\sum_{t\geq 0}{t\rank M_x(t)}.
$$
We let $\bb_x(M)$ denote the largest break, i.e., the largest $t\geq
0$ such that $M_x(t)\not=0$. Notice then the trivial inequalities 
\begin{align}
\swan_x(M)&\leq (\rank M)\bb_x(M),
\label{eq-sw-up}
\\
\bb_x(M)&\leq \swan_x(M).
\label{eq-sw-do}
\end{align}
\par
In addition, if $M=M_1\otimes M_2$ we have~\cite[Lemma
1.3]{katz-gsksm} 
\begin{equation}\label{eq-sw-tenseur}
\bb_x(M)\leq \max(\bb_x(M_1),\bb_x(M_2))\leq \bb_x(M_1)+\bb_x(M_2).
\end{equation}
\par
If $M$ is a finite dimensional $\Qq_{\lambda}$-vector space, with
$\Qq_{\lambda}$ a finite extension of $\Qq_{\ell}$ with ring of integers
$\Zz_{\lambda}$, for some $\ell\not=p$, and if
$\mathbf{M}\subset M$ is an invariant $\Zz_{\lambda}$-lattice with reduction
$\mathbf{M}/\lambda\mathbf{M}$, then we have~\cite[Rem. 1.10]{katz-gsksm}
\begin{equation}\label{eq-eq-swans}
\swan_x(M)=\swan_x(\mathbf{M})=\swan_x(\mathbf{M}/\ell\mathbf{M}).
\end{equation}
\par
Finally, the main reason the Swan conductor enters in our computation
is the fundamental formula of Grothendieck-Ogg-Shafarevitch:
\begin{proposition}
If $\Qq_{\lambda}$ is a finite extension of $\Qq_{\ell}$ and
$\sheaf{F}$ is a lisse $\Qq_{\lambda}$-sheaf on $U$ of rank $r$, we have
\begin{equation}\label{eq-e-p}
\chi_c(U,\sheaf{F})=
r\chi_c(U,\Qq_{\ell})-\sum_{x\in S}{\swan_x(\sheaf{F})},
\end{equation}
where $\swan_x(\sheaf{F})$ is $\swan_x(M)$ for the $\Qq_{\lambda}$-vector
space which is the representation space of the representation
corresponding to $\sheaf{F}$.
\end{proposition}
See e.g.~\cite[2.3.1, 2.3.3]{katz-gsksm} for a sketch of the proof.

\begin{proof}[Proof of Proposition~\ref{pr-boundcurve}]
Since $U$ is affine and smooth we have
$h^0_c(U,\pi(\rho))=0$ and
$$
\sigma'_c(U,\pi(\rho))=h^1_c(U,\pi(\rho)),
$$
while the Euler-Poincar\'e characteristic is 
$$
\chi_c(U,\pi(\rho))=-h^1_c(U,\pi(\rho))+
h^2_c(U,\pi(\rho)).
$$
\par
We want to bound $-\chi_c(U,\pi(\rho))$, and for this start from the
Euler-Poincar\'e formula~(\ref{eq-e-p}) for $\pi(\rho)$ (which takes
value in some $GL(r,\Qq_{\lambda})$):
$$
\chi_c(U,\pi(\rho))=(\dim \pi)\chi_c(U,\Qq_{\ell})-
\sum_{x\in S}{\swan_x(\pi(\rho))}.
$$
By~(\ref{eq-sw-up}) we get bounds
$$
\swan_x(\pi(\rho))\leq (\dim \pi)\bb_x(\pi\circ\rho)
\leq (\dim \pi)\bb_x(M),
$$
where $M$ is the $\Zz[1/p]$-module 
$$
M=M_1\otimes \cdots\otimes M_k,
$$
with $M_i\simeq \Fp_{\lambda_i}^r$, with the action
$\rho=\rho_1\otimes \cdots\otimes\rho_k$ of $\pi_1(U,\barre{\eta})$.
The last inequality is simply because the action of the
inertia group ``on'' $\pi\circ\rho$ factors through that on $M$.
\par
Now we have by~(\ref{eq-sw-tenseur}) and~(\ref{eq-sw-do})
$$
\bb_x(M)\leq \max\bb_x(\rho_i)\leq \sum_{i}{\bb_x(\rho_i)}\leq
\sum_i{\swan_x(\rho_i)}.
$$
Hence
$$
\swan_x(\pi(\rho))\leq (\dim \pi)\sum_{i}{\swan_x(\rho_i)}.
$$
Now we can use the fact that each $\rho_i$ is the reduction of the
$\Zz_{\lambda_i}$-sheaf $\sheaf{F}_i$ and use~(\ref{eq-eq-swans}) to get
$$
\sum_{i}{\swan_x(\rho_i)}=\sum_{i}{\swan_x(\sheaf{F}_i\otimes
\Qq_{\lambda_i})},
$$
and by~(\ref{eq-e-p}) again we have for all $i$
$$
\sum_{x\in S}{\swan_x(\sheaf{F}_i\otimes\Qq_{\lambda_i})}=
r\chi_c(U,\Qq_{\lambda_i})-\chi_c(U,\sheaf{F}_i).
$$
The crucial point is that $\chi_c(U,\sheaf{F}_i)$ is independent
of $i$ because the sheaves $\sheaf{F}_i$ form a compatible system
(it is minus the degree of the common $L$-function of the sheaves
$\sheaf{F}_{i}$), and so of course is $\chi_c(U,\Qq_{\lambda_i})$. So
the sum of the Swan conductors of the 
$\sheaf{F}_i\otimes\Qq_{\lambda_i}$ is independent of $i$. We denote
this common value by $w$, and thus we get
\begin{align*}
-\chi_c(U,\pi(\rho))&\leq 
(\dim\pi)\Bigl\{
-\chi_c(U,\Qq_{\ell})
+\sum_x{\sum_i{\swan_x(\sheaf{F}_i\otimes\Qq_{\lambda_i})}}
\Bigr\}\\
&= (\dim\pi)\Bigl\{
-\chi_c(U,\Qq_{\ell})
+kw\Bigr\}
%
\end{align*}
\par
We add the requisite $h^2_c(U,\pi(\rho))$ which is
trivially $\leq \dim\pi$ by the co-invariant description
$$
H^2_c(U,\sheaf{F})\simeq \sheaf{F}_{\pi_1(U,\barre{\eta})}(-1)
$$
for any lisse $\barre{\Qq}_{\ell}$-sheaf $\sheaf{F}$, and therefore we get
\begin{equation}\label{eq-w}
\sigma'_c(U,\pi(\rho))=h^2_c(U,\pi(\rho))
-\chi_c(U,\pi(\rho))
\leq (\dim \pi)(1+C)\text{ with }
C=-\chi_c(U,\Qq_{\ell})+kw.
\end{equation}
\end{proof}

We now come to the result that will be used for the case where
Assumption (ii) of Theorem~\ref{th-main} holds.
We will use the following result of Katz, building on work of Bombieri
and Adolphson and Sperber:

\begin{proposition}\label{pr-katz}
Let $q=p^k$, $U/\barre{\Fp}_q$ a
smooth affine connected scheme of dimension $d\geq 1$ which can be
embedded in $\Aa^N$ as a closed subscheme defined by the vanishing of
$r$ polynomials of degree $\leq \delta$. Then we have
$$
\sigma_c(U,\Qq_{\ell})\leq A(N,r,\delta)
$$
for some constant $A(N,r,\delta)$;
one can take 
\begin{equation}\label{eq-arj}
A(N,r,\delta)=2^r6(3+r\delta)^{N+1}.
\end{equation}
\end{proposition}

This is Theorem~2 of~\cite{katz-betti} together with its corollaries.

\begin{proposition}\label{pr-sigmac}
Let $q=p^k$, $U/\barre{\Fp}_q$ a
smooth affine connected scheme over $\barre{\Fp}_q$ of dimension $d\geq 1$,
$\varphi\,:\,V\ra U$ a finite étale connected Galois covering of
degree prime to $p$. There exists a constant $C(U)$ such that 
\begin{equation}\label{eq-sigmac}
\sigma_c(V,\Qq_{\ell})\leq C(U)(\deg \varphi).
\end{equation}
More precisely, if $d=1$ one can take $C(U)=\sigma_c(U,\Qq_{\ell})$. If 
$d\geq 2$, let $N$, $r$, $\delta$ be as in Proposition~\ref{pr-katz}
for $U$. Then one can take $C(U)=C(N,r,\delta)$, 
where
\begin{equation}\label{eq-rec}
C(N,r,\delta)=2\sum_{j=1}^{N-1}{A(j,r,\delta)}+A(N,r,\delta)\leq
12N2^r(3+r\delta)^{N+1}.
\end{equation}
\end{proposition}

Something like this may be already known but we haven't found it in the
literature. 
The proof will proceed by induction on $d$, following the method used
by Katz in~\cite[Th. 2]{katz-betti}. For the induction step we need
the following version of an affine Lefschetz theorem:

\begin{proposition}\label{pr-lefschetz}
Let $U/\barre{\Fp}_q$ be a
smooth connected affine scheme over $\barre{\Fp}_q$ of dimension $d\geq 2$,
$\varphi\,:\,V\ra U$ a 
finite étale connected Galois covering with Galois group $G$. Fixing an
immersion $i\,:\, U\ra \Aa^N$ for some $N\geq 1$, there exists an
affine hyperplane $H\subset \Aa^N$ such that $U\cap H$ is connected
and smooth, $W=\varphi^{-1}(U\cap H)$ is connected and smooth, and in
the diagram  
$$
\begin{CD}
  W @>>> V\\
@V\varphi_1 VV  @VV\varphi V \\
U\cap H @>>> U
\end{CD}
$$
the map $\varphi_1$ is a finite étale Galois covering with group $G$
and the induced maps in étale cohomology
\begin{align}
H^i(V,\Qq_{\ell})&\ra H^{i}(W,\Qq_{\ell})\label{eq-2}
\end{align}
are isomorphisms for $i<d-1$ and injective for $i=d-1$.
\end{proposition}

\begin{proof}
For any hyperplane $H$, it is of course true that $W\ra U\cap H$ is a
finite étale covering with Galois group $G$, possibly
disconnected. However there exists an open dense set of hyperplanes
$H$ for which $W$ is indeed connected by~\cite[Cor. 3.4.2]{katz-act}
with the data $(k,E,f,\pi)=(\barre{\Fp}_q,V,0,i)$.
\par
Further, the existence of an open dense set of hyperplanes $H$ such that the
induced maps 
$$
H^i(V,\Qq_{\ell})\ra H^{i}(\varphi^{-1}(U\cap H),\Qq_{\ell})
=H^i(W,\Qq_{\ell})
$$
satisfy the required condition is the special case
of~\cite[Cor. 3.4.1]{katz-act} for the data
$(k,E,f,\pi)=(\barre{\Fp}_q,V,0,i\circ\varphi)$ (compare the proof 
of~\cite[Cor. 3.4.2]{katz-act}). The existence of a third open dense
set of $H$ for which $U\cap H$ is smooth connected
is Cor. 3.4.3 of loc. cit. 
\par
Intersecting those three open
dense subsets of hyperplanes, one finds one where all the required
conditions hold. 
\end{proof}

\begin{proof}[Proof of Proposition~\ref{pr-sigmac}]
First because $U$ is smooth affine and $\varphi$
étale, hence finite, $V$ is also affine and smooth. 
\par
We recall now some deep facts about étale cohomology. First,
since $V$ is smooth and connected we have
$\sigma_c(V,\Qq_{\ell})=\sigma(V,\Qq_{\ell})$ by Poincar\'e
duality~(see e.g.~\cite[VI.3]{sga4half}). 
\par
Next, by the affine cohomological dimension theorem we
have 
\begin{equation}\label{eq-aff-dim}
H^i(U,\Qq_{\ell})=H^i(V,\Qq_{\ell})=0\text{ for } i>d,
\end{equation}
see e.g.~\cite[IV.6.4]{sga4half}.
\par
Finally, because $\varphi$ is an étale Galois covering of degree prime
to $p$, it is moderately ramified and we have 
\begin{equation}\label{eq-chi-et}
\chi(V,\Qq_{\ell})=(\deg \varphi)\chi(U,\Qq_{\ell})
\end{equation}
which is due to
Deligne-Lusztig for $\chi_c$ (see~\cite[2.6, Cor. 2.8]{illusie}), and
we have $\chi=\chi_c$ for $U$ and $V$ as proved by
Laumon~\cite{laumon}.
\par
Now we are ready to start the proof by induction. Consider first
$d=1$. We have by~(\ref{eq-aff-dim}) and Poincaré duality which gives
$h^0(U,\Qq_{\ell})=h^0(V,\Qq_{\ell})=1$ that
$$
\sigma(U,\Qq_{\ell})=2-\chi(U,\Qq_{\ell}),\quad
\sigma(V,\Qq_{\ell})=2-\chi(V,\Qq_{\ell}).
$$
By~(\ref{eq-chi-et}) we get
$$
\sigma(V,\Qq_{\ell})=2-(\deg\varphi)\chi(U,\Qq_{\ell})
\leq
(\deg\varphi)(2-\chi(U,\Qq_{\ell}))=\deg(\varphi)\sigma(U,\Qq_{\ell})
$$
so that we can indeed take $C(U)=\sigma_c(U,\Qq_{\ell})$ in that
case. This is the (first) conclusion required for $d=1$. The
alternative bound is also valid since $\sigma(U,\Qq_{\ell})\leq
A(N,r,\delta)\leq C(N,r,\delta)$ with $N$, $r$, $\delta$ as described,
by Proposition~\ref{pr-katz} and $C(N,r,\delta)$ defined
by~(\ref{eq-rec}). 
\par
Now assume that $\dim U=d$ and~(\ref{eq-sigmac}) holds for dimension
$d-1$ with the constant~(\ref{eq-rec}). Fix an embedding $i\,:\, U\ra
\Aa^N$ for some $N$ (with the 
attending $r$ and $\delta$). By Proposition~\ref{pr-lefschetz} there
exists an hyperplane $H\subset \Aa^N$ 
such that the maps~(\ref{eq-2}) are, in particular,
all injective for $i\leq d-1$. This implies by~(\ref{eq-aff-dim})
$$
\sigma(V,\Qq_{\ell})\leq \sigma(W,\Qq_{\ell})+h^d(V,\Qq_{\ell})
$$
on the one hand, and on the other hand we find
$$
h^d(V,\Qq_{\ell})\leq h^d(V,\Qq_{\ell})+
h^{d-1}(W,\Qq_{\ell})-h^{d-1}(V,\Qq_{\ell})
=(-1)^d\chi(V,\Qq_{\ell})+(-1)^{d-1}\chi(W,\Qq_{\ell}),
$$
so that altogether we have the inequality
$$
\sigma(V,\Qq_{\ell})\leq (-1)^d\chi(V,\Qq_{\ell})+
(-1)^{d-1}\chi(W,\Qq_{\ell})+\sigma(W,\Qq_{\ell}).
$$
\par
Now using twice~(\ref{eq-chi-et}) for $V\ra U$ and $W\ra U\cap H$,
which are both Galois with group $G$ of order prime to $p$, we find
$$
\sigma(V,\Qq_{\ell})\leq (\deg\varphi)((-1)^d\chi(U,\Qq_{\ell})+
(-1)^{d-1}\chi(U\cap H,\Qq_{\ell}))+
\sigma(W,\Qq_{\ell}).
$$
By induction applied to $W\ra U\cap H$, since $U\cap H$ is embedded in
$H\simeq \Aa^{N-1}$ using the same number and degree of polynomials as
$U$, we can estimate the last term by 
$$
\sigma(W,\Qq_{\ell})\leq (\deg \varphi)C(N-1,r,\delta)
$$
and get 
$$
\sigma(V,\Qq_{\ell})\leq (\deg\varphi)\Bigl\{(-1)^d\chi(U,\Qq_{\ell})+
(-1)^{d-1}\chi(U\cap H,\Qq_{\ell})+
C(N-1,r,\delta)\Bigr\}.
$$
Since 
\begin{gather*}
|\chi(U,\Qq_{\ell})|\leq \sigma(U,\Qq_{\ell})\leq A(N,r,\delta)
\\
|\chi(U\cap H,\Qq_{\ell})|\leq \sigma(U\cap H,\Qq_{\ell})
\leq A(N-1,r,\delta),
\end{gather*}
by Proposition~\ref{pr-katz}, we get
$$
\sigma(V,\Qq_{\ell})\leq (\deg\varphi)\Bigl\{A(N,r,\delta)+
A(N-1,r,\delta)+C(N-1,r,\delta)\Bigr\}=
(\deg\varphi)C(N,r,\delta),
$$
which is the result for $U$. The last estimate in~(\ref{eq-rec}) is a
crude consequence of the corresponding one for $A(j,r,\delta)$ given
in~(\ref{eq-arj}). 
\end{proof}

\begin{proposition}\label{pr-boundok}
Let $U/\barre{\Fp}_q$ be a smooth affine connected scheme of
dimension $d\geq 1$. Let
$\rho\,:\, \pi_1(U,\barre{\eta})\ra G$ be a surjective homomorphism
with $G$ finite of order prime to $p$, let $\pi\,:\, G\ra
GL(r,\barre{\Qq}_{\ell})$ be a representation of $G$ and
$\pi(\rho)=\pi\circ\rho$ the corresponding lisse sheaf on $U$.
There exists a constant $C(U)$ depending only on $U$ such that
\begin{equation}\label{eq-boundok}
\sigma_c(U,\pi(\rho))\leq C(U)|G|(\dim\pi).
\end{equation}
\end{proposition}

\begin{proof}
(Compare with (4) of Theorem~9.2.6 of Katz and
Sarnak~\cite{katz-sarnak})
Let $\varphi\,:\,V\ra U$ be the connected étale covering with group
$G$ corresponding to the kernel of $\rho$. It follows that
$\varphi^*(\pi(\rho))$ is trivial on $V$, i.e., seeing $\pi(\rho)$ as a
$\Qq_{\lambda}$-lisse sheaf, where $\Qq_{\lambda}$
is a finite extension of $\Qq_{\ell}$ for which $\pi$ has image in
$GL(r,\Qq_{\lambda})$, we have
$$
\varphi^*\pi(\rho)\simeq \Qq_{\lambda}^{r}.
$$
\par
Since $V\ra U$ is étale and Galois, the Galois group $G$ acts on the
cohomology groups $H^i_c(V,\varphi^*\pi(\rho))$ and we have (by the
Hochschild-Serre spectral sequence for $V\ra U$ for instance) for all
$i$  
$$
H^i_c(V,\varphi^*\pi(\rho))^G
\simeq H^i_c(U,\pi(\rho))
$$
hence 
$$
\sigma_c(U,\pi(\rho))\leq \sigma_c(V,\varphi^*\pi(\rho))=
\sigma_c(V,\Qq_{\lambda}^{r})=
r\sigma_c(V,\Qq_{\lambda})=r\sigma_c(V,\Qq_{\ell})
$$
by the formula $H^i_c(V,\Qq_{\lambda})=H^i_c(V,\Qq_{\ell})\otimes
\Qq_{\lambda}$. 
\par
Since the group $G$ is assumed to have order prime to $p$, we have by
Proposition~\ref{pr-sigmac} 
$$
\sigma_c(V,\Qq_{\ell})\leq C(U)|G|
$$
for some constant $C(U)$ independent of $\pi$, and the proposition
follows by combining these two inequalities.
\end{proof}

\begin{remark}
Contrary to our first optimistic version, the
condition that $(\deg\varphi,p)=1$ in Proposition~\ref{pr-sigmac} is
certainly necessary, as the following example (communicated by Katz)
shows: take $U$ to be the affine line $\Aa^1 $ 
with coordinate $x$, and take $V=V_d$ to be the curve $y^p - y = x^d$,
for $d$ prime to $p$. Then we have $\sigma_c(V_d)=1 + (p-1)(d-1)$. So as $d$
grows, although the degree of the covering $V_d\ra U$ stays $p$, we
see that $\sigma_c(V_d)$ is unbounded. 
\par
Since the covering is 
also Galois with group $\Zz/p\Zz$ in this case, this also shows that
Proposition~\ref{pr-boundok} does not extend to arbitrary groups: the
covering $V\ra U$ corresponds to a surjective map $\rho\,:\, \pi_1(U)\ra
\Zz/p\Zz$, the representations $\pi=\psi$ are the additive characters
of $\Zz/p\Zz$, and we have
$$
\sigma_c(V_d,\Qq_{\ell})=\sum_{\psi}{\sigma_c(U,\psi(\rho))}
$$
(which amounts to the standard counting of points on
$V_d(\Fp_{q^n})$ by means of additive character sums, or the
construction of the sheaves corresponding to those sums) and therefore
a bound like~(\ref{eq-boundok}) would give $\sigma_c(V_d)\leq
p^2C(U)$, which is also incorrect.
\par
The condition that the covering be Galois is necessary for the proof
of Proposition~\ref{pr-sigmac} because otherwise~(\ref{eq-chi-et})
may fail, even for a covering of degree prime to $p$, as in the following
example (again communicated by Katz): take  the finite étale covering
$\mathbf{G}_m\ra \Aa^1$ (over $\Fp_q$) given by $x\mapsto x^p + 1/x$. We
have $\chi(\Aa^1)=1$, whereas $\chi(\mathbf{G}_m)=0$, and the covering
is of degree $p+1$.
\par
Still one may hope that an analogue of Proposition~\ref{pr-boundcurve}
holds for arbitrary $U$, which would give a corresponding general
version of Theorem~\ref{th-main} and its applications.
\end{remark}

\section{Proof of the bilinear form estimate}
\label{sec-bilinear}

We come back to the notation of Section~\ref{sec-main} before and in
the statement of Theorem~\ref{th-main}, which we will now prove.
\par
The analytic principle for the proof of Theorem~\ref{th-main} is quite
simple and very well established in analytic number theory. We proceed
by duality, as first conceived by Vinogradov: for given $L\geq 1$ and
$\Delta\geq 0$, it is equivalent to prove that
$$
\sum_{\ell\leq L}{\sums_{\pi\not=1}{
\Bigl|
\sum_{u\in U(\Fp_q)}{\alpha(u)\Tr(\pi\circ\rho_{\ell})(\frob_u)}
\Bigr|^2
}}
\leq \Delta \sum_{u\in U(\Fp_q)}{|\alpha(u)|^2},
$$
for arbitrary $\alpha(u)\in\Cc$, or to prove that
\begin{equation}\label{eq-dual}
\sum_{u\in U(\Fp_q)}{
\Bigl|
\sum_{\ell\leq L}{\sums_{\pi\not=1}{\beta(\ell,\pi)
\Tr(\pi\circ\rho_{\ell})(\frob_u)}}
\Bigr|^2}
\leq \Delta \sum_{\ell\leq L}{\sums_{\pi\not=1}{|\beta(\ell,\pi)|^2}}.
\end{equation}
for arbitrary $\beta(\ell,\pi)\in\Cc$. Recall that $\pi$
runs over a set $\Pi_{\ell}$ of irreducible representations of
$G_{\ell}$ up to twist by characters of $\Gamma_{\ell}$.
\par
The dual form is more manageable here. Denote by $\ideal{B}(\beta)$
the left-hand side of~(\ref{eq-dual}). Expanding the 
square we get
\begin{equation}\label{eq-expand}
\ideal{B}(\beta)=\sum_{\ell\leq L}{\sums_{\pi\not=1}{\sum_{\ell'\leq L}{
\sums_{\pi'\not=1}{
\beta(\ell,\pi)\overline{\beta(\ell',\pi')}
\mathfrak{S}(\ell,\pi;\ell',\pi')}}}},
\end{equation}
with
$$
\mathfrak{S}(\ell,\pi;\ell',\pi')=
\sum_{u\in U(\Fp_q)}{
\Tr(\pi\circ\rho_{\ell})(\frob_u)
\overline{
\Tr(\pi'\circ\rho_{\ell'})(\frob_u)
}}.
$$
The crucial point is the following estimation for the individual
$\mathfrak{S}(\ell,\pi;\ell',\pi')$. 

\begin{proposition}\label{pr-coeff}
With notation as above, and in particular under the assumption that
the sheaves are \disj .
\par
\emph{(i)} If the monodromy groups $G^g_{\ell}$ are of
prime-to-p order, we have
\begin{align*}
|\mathfrak{S}(\ell,\pi;\ell,\pi)-|\hat{\Gamma}_{\ell}^{\pi}|
q^d|&
\leq q^{d-1/2}|G_{\ell}|(\dim \pi)^2C(\barre{U}),\\ 
|\mathfrak{S}(\ell,\pi;\ell,\pi')|&
\leq q^{d-1/2}|G_{\ell}|(\dim \pi)
(\dim \pi')C(\barre{U}),\text{ if }\pi\not=\pi'\\
|\mathfrak{S}(\ell,\pi;\ell',\pi')|&
\leq q^{d-1/2}|G_{\ell}||G_{\ell'}|
(\dim \pi)(\dim \pi')C(\barre{U}),
\text{ if }\ell\not=\ell',
\end{align*}
where $C(\barre{U})$ is given by Proposition~\ref{pr-sigmac}.
\par
\emph{(ii)} If $U$ is a curve, and the sheaves arise from a compatible
system of $\Zz_{\lambda}$-sheaves $\sheaf{F}_{\ell}$, we have
\begin{align*}
|\mathfrak{S}(\ell,\pi;\ell,\pi)-|\hat{\Gamma}_{\ell}^{\pi}|
q^d|&
\leq q^{d-1/2}(\dim \pi)^2D(\barre{U},(\sheaf{F}_{\ell})),\\ 
|\mathfrak{S}(\ell,\pi;\ell',\pi')|&
\leq q^{d-1/2}(\dim \pi)(\dim \pi')D(\barre{U},(\sheaf{F}_{\ell})),
\text{ if }\ell\not=\ell'\text{ or }\pi\not=\pi',
\end{align*}
where $D(\barre{U},(\sheaf{F}_{\ell}))$ is the constant
$C(\barre{U},(\sheaf{F}_{\ell}),2)$ of Proposition~\ref{pr-boundcurve}.
\end{proposition}

Taking this for granted, we finish quickly the proof of
Theorem~\ref{th-main}. By~(\ref{eq-expand}) we have
trivially~(\ref{eq-dual}) with 
$$
\Delta=\max_{\ell,\pi}{
\sum_{\ell'}{\sums_{\pi'\not=1}{|\mathfrak{S}(\ell,\pi;\ell',\pi')|}}
},
$$
and by Proposition~\ref{pr-coeff}, we thus get~(\ref{eq-dual}) with
\begin{equation}\label{eq-delta-final}
\Delta=\max_{\ell,\pi}\Bigl\{|\hat{\Gamma}_{\ell}^{\pi}|
q^{d}+q^{d-1/2}C(\barre{U})
|G^g_{\ell}|(\dim \pi)\Bigl\{
\sums_{\pi'}{(\dim \pi')}+
\sum_{\ell'\not=\ell}{|G^g_{\ell'}|\sums_{\pi'\not=1}{(\dim \pi')}}
\Bigr\}
\Bigr\}
\end{equation}
in the case of monodromy of order prime to $p$, and
\begin{equation}\label{eq-delta-curve}
\Delta=\max_{\ell,\pi}\Bigl\{
|\hat{\Gamma}_{\ell}^{\pi}|
q^{d}+q^{d-1/2}D(\barre{U},(\sheaf{F}_{\ell}))
(\dim \pi)\Bigl\{
\sums_{\pi'}{(\dim \pi')}+
\sum_{\ell'\not=\ell}{\sums_{\pi'\not=1}{(\dim \pi')}}
\Bigr\}
\Bigr\}
\end{equation}
in the case of a curve with a compatible system. We estimate
all those terms in terms of the parameters $s$ and $t$ of~(\ref{eq-s})
using Lemma~\ref{lm-sizes},~(1).
In the first case we obtain by appealing also
to~(\ref{eq-kappa}) and to~(\ref{eq-s}) that
\begin{equation}\label{eq-delta-final-2}
\Delta\leq \kappa q^d+2q^{d-1/2}C(\barre{U})(c_1^5c_2)^{1/2}L^{1+5s/2+t/2}.
\end{equation}
\par
In the second case we obtain similarly
\begin{equation}\label{eq-delta-curve-2}
\Delta\leq \kappa q^d+2q^{d-1/2}D(\barre{U},(\sheaf{F}_{\ell}))
c_1c_2^{1/2}L^{1+s+t/2}.
\end{equation}
Thus Theorem~\ref{th-main} follows by duality.

\begin{proof}[Proof of Proposition~\ref{pr-coeff}]
The proof is now an easy application of the Grothendieck-Lefschetz
trace formula and Deligne's main theorem of~\cite{weil2} (compare with
\cite[p. 162,163]{chavdarov}).
The only subtlety is that the dependency of the error terms on
$\ell$, $\ell'$, $\pi$, $\pi'$, must remain controlled, 
and for this we need the results of Section~\ref{sec-betti}.
\par
If $\ell=\ell'$ we let $G=G_{\ell}$, $G^g=G^g_{\ell}$. The
representation $\rho_{\ell}$ gives a surjective map
$\pi_1(U,\barre{\eta})\ra G$. Let
$\tau=\pi\otimes\tilde{\pi}'$. This is a (not necessarily
irreducible) representation of $G$, and we will consider the sheaf 
$\sheaf{F}=\tau\circ \rho_{\ell}$, which is of the type considered
in Proposition~\ref{pr-boundok} and Proposition~\ref{pr-boundcurve}
(after seeing the representation $\tau$ as taking value in
$GL(r,\barre{\Qq}_{\ell})$, as we can since it is a representation in
characteristic $0$).
\par
If $\ell\not=\ell'$, we let $G=G_{\ell}\times G_{\ell'}$,
$G^g=G^g_{\ell}\times G^g_{\ell'}$. By the assumption that the family
of sheaves is \disj , the product map
$(\rho_{\ell},\rho'_{\ell})$ is still a surjective map
$$
\pi_1(\barre{U},\barre{\eta})\fleche{(\rho_{\ell},\rho_{\ell'})} G^g.
$$
Let $\tau(g,g')=\pi(g)\otimes\tilde{\pi}'(g')$ (the ``external''
product), so $\tau$ is an irreducible representation of $G$. We will 
consider the sheaf $\sheaf{F}=\tau\circ (\rho_{\ell},\rho'_{\ell})$,
again of the type considered in Proposition~\ref{pr-boundok} and
Proposition~\ref{pr-boundcurve}. 
\par
In both cases, because $G$ is a finite group, the eigenvalues of the
image of $\tau$ are roots of unity so $\sheaf{F}$ is punctually pure
of weight $0$.
\par
Also in either case,
the main point is that the local trace at
$u\in U(\Fp_q)$ of $\sheaf{F}$ is given by construction by
$$
\Tr(\frob_u\,\mid\, \sheaf{F})=
\Tr(\pi\circ\rho_{\ell})(\frob_u)
\overline{\Tr(\pi'\circ\rho_{\ell'})(\frob_u)}
$$
and therefore the fundamental Grothendieck-Lefschetz Trace Formula
(see~\cite{grothendiec},~\cite{sga4half},~\cite[VI.13]{milne}) states
that
$$
\mathfrak{S}(\ell,\pi;\ell',\pi')=
\sum_{u\in U(\Fp_q)}{\Tr(\frob_u\,\mid\,\sheaf{F})}
=\sum_{0\leq i\leq 2d}
{(-1)^i\Tr (\frob\,\mid\, H^i_c(\barre{U},\sheaf{F}))}.
$$
\par
By Deligne's Weil II Theorem~\cite[p. 138]{weil2}, the eigenvalues of
the geometric Frobenius automorphism $\frob$ acting on
$H^i_c(\barre{U},\sheaf{F})$ are algebraic 
integers all conjugates of which are of absolute value $\leq
q^{i/2}$. 
\par
It is easy to compute $H^{2d}_c(\barre{U},\sheaf{F})$ (the action on
which contributes potentially terms of maximal size $q^d$),
using the formula
$$
H^{2d}_c(\barre{U},\sheaf{F})=V_{\pi_1(\barre{U},\barre{\eta})}(-d)=
W_{G^g}(-d)
$$
where $V=\sheaf{F}_{\barre{\eta}}$ is the space on which the
representation which ``is'' 
$\sheaf{F}$ acts and $W$ is the space of the representation $\tau$ of
$G^g$. The second equality above holds because the
disjointness condition shows that the map
$\pi_1(\barre{U},\barre{\eta}) \ra G^g$ through which the action factors
is always surjective (as already observed previously).
\par
The crucial point is that this coinvariant space is zero unless
$\ell=\ell'$ and $\pi=\pi'$. Indeed, if $\ell\not=\ell'$, decomposing
$\pi$ and $\pi'$ restricted to $G_{\ell}^g$ and $G_{\ell'}^g$ as sums
of irreducible representations, the dimension of the coinvariant space
(which is the same as that of the invariants under $G^g$ because we are
working with finite groups) is the sum of the dimension of invariants
for the pairwise tensor products of the components; but those are
non-trivial irreducible representations of $G^g=G_{\ell}^g\times
G_{\ell'}^g$ so each term of the sum is zero.
\par
If $\ell=\ell'$, the last statement in Lemma~\ref{lm-ortho} exactly
says that the space of invariants is of dimension
$|\hat{\Gamma}_{\ell}^{\pi}|$ for $\pi=\pi'$ and $0$ otherwise; this
is where it is necessary to restrict to representations unrelated by
twists.
\par
Thus we derive the bound
$$
\Bigl|\sum_{u\in U(\Fp_q)}{\Tr(\pi\circ\rho_{\ell})(\frob_u)
\overline{\Tr(\pi'\circ\rho_{\ell'})(\frob_u)}}\Bigr|
\leq q^{d-1/2}\sigma'_c(\barre{U},\sheaf{F})
$$
if $(\ell,\pi)\not=(\ell',\pi')$ and
$$
\Bigl|\sum_{u\in U(\Fp_q)}{\Tr(\pi\circ\rho_{\ell})(\frob_u)
\overline{\Tr(\pi'\circ\rho_{\ell'})(\frob_u)}}-
|\hat{\Gamma}_{\ell}^{\pi}|q^d\Bigr|
\leq q^{d-1/2}\sigma'_c(\barre{U},\sheaf{F})
$$
otherwise.
\par
Inserting the bounds for $\sigma'_c(\barre{U},\sheaf{F})$ from
Proposition~\ref{pr-boundcurve} or Proposition~\ref{pr-boundok}
respectively, and looking at the various cases, the
proposition follows.  
\end{proof}

\begin{remark}
If the conditions of Theorem~\ref{th-main} are not satisfied, we see
that we still get an inequality
$$
\sum_{\ell}{\sums_{\pi}{\Bigl|\sum_{u\in U(\Fp_q)}{
\alpha(u) \Tr(\pi\circ \rho_{\ell})(\frob_u)}\Bigr|^2}}
\leq (\kappa q^d+q^{d-1/2}D)\sum_{u\in U(\Fp_q)}{|\alpha(u)|^2}
$$
for fixed $L$, with
$$
D=\max_{\ell,\pi}{\sum_{\ell'}{\sums_{\pi'\not=1}{\sigma'_c(\barre{U},
\sheaf{F}_{\tau_{\pi,\pi'}})}}}.
$$
The point is that $D$ is independent of $q$ so this is still
non-trivial when applied for $U\times\Fp_{q^n}$ with $n\ra +\infty$.
In a large sieve context as in Proposition~\ref{pr-large-sieve}, it
leads to
$$
\limsup_{n\ra +\infty}{\frac{|S(U\times \Fp_{q^n},\Omega;L)|}{q^{nd}}}
\leq \frac{1}{P(L)}
$$
for fixed $L$, and using the same sieve as we will in the proof of
Theorem~\ref{th-chavdarov} for all $L\geq 2$ and taking $L\ra
+\infty$, this recovers Chavdarov's irreducibility
theorem~\cite[Th. 2.3]{chavdarov} for an 
arbitrary family $C/U$ of genus $g$ curves with geometric monodromy
modulo $\ell$ equal to $Sp(2g,\Fp_{\ell})$ for almost all $\ell$:
$$
\lim_{n\ra +\infty}{\frac{
|\{
u\in U(\Fp_{q^n})\,\mid\,
\det(1-T\frob_u\,\mid\, H^1(\barre{C}_u,\Zz_{\ell}))
\text{ has small Galois group}
\}|
}{|U(\Fp_{q^n})|}}=0.
$$
\end{remark}

\section{Zeta functions of families of curves}\label{sec-chavdarov}

We now come to the application of the large sieve to a strong form of
Chavdarov's Theorem on the generic behavior of the numerators of zeta
functions of curves in families. If the genus is fixed, most of the
work is already done in the previous sections or in Chavdarov's paper,
but we will look for arguments uniform with respect to $g$ so that, in
some cases at least, we obtain results valid even for $g$ large
(though not for $q$ fixed, $g\ra +\infty$).
\par
First, we recall the definition of the zeta function of a curve over a
finite field, in concrete terms (so the statements at least can be
understood without knowledge of étale cohomology), by recalling the
diophantine meaning of the polynomials involved.  
\par
Let $C/\Fp_q$ be a smooth projective curve of genus $g$ over a finite
field (all curves here and below are assumed to be geometrically
connected). Its zeta function $Z(C)$ is the formal 
power series given by the diophantine definition 
$$
Z(C)=\exp\Bigl(\sum_{n\geq 1}{\frac{|C(\Fp_{q^n})|}{n}T^n}\Bigr),
$$
where $|C(\Fp_{q^n})|$ is the number of ``solutions'' to the equations
which define $C$ with coordinates
in the extension field $\Fp_{q^n}$. A fundamental result due to
F.K. Schmidt in this 
case is that there exists a polynomial $P_C\in\Zz[T]$ of degree $2g$
with $P_C(0)=1$ such that
$$
Z(C)=\frac{P_C(T)}{(1-T)(1-qT)}.
$$
The cohomological definition is that the polynomial $P_C(T)$ can
be described as the (reversed) characteristic polynomial of the
geometric Frobenius automorphism acting on a suitable étale
cohomology group, specifically 
\begin{equation}\label{eq-cohomo-interp}
P_C(T)=\det(1-T\frob\,\mid\, H^1(\barre{C},\Zz_{\ell})).
\end{equation}
The question investigated by Chavdarov concerns the splitting field of
this integer polynomial as $C$ varies in an algebraic family, e.g. in
a hyperelliptic family
$$
C_u\,:\,y^2=f(x)(x-u)
$$
where $f$ is a fixed polynomial in $\Fp_q[X]$ of degree $2g$ with
distinct roots in $\barre{\Fp}_q$, and $u$ is the parameter that can take
any value in $\barre{\Fp}_q$ which is not a zero of $f$ (these
conditions ensure that the curve $C_u$ suitably ``compactified'' is a
smooth projective curve of the given genus $g\geq 1$).
\par
There is an a-priori condition on the splitting field of
the polynomial $P_C$ because it satisfies the ``functional equation''
$$
T^{2g}P_C\Bigl(\frac{q}{T}\Bigr)=P_C(T),
$$
(or equivalently, if $\alpha\in\Cc$ is a root of $P_C$, then
$q\alpha^{-1}$ is also a root). This means
that the ``splitting algebra''
$\Qq[T]/(f)$ has Galois group $G$ which can be seen as a subgroup of the
group $W_{2g}$ of signed permutations of $\{1,\ldots,2g\}$.
In other words, $W_{2g}$ is the group of permutations of $g$ pairs of
elements preserving the pairs. In particular, if the polynomial is
irreducible, its splitting field has maximal Galois group $G\simeq
W_{2g}$ if and only if the splitting field is of maximal degree
$|W_{2g}|=2^g g!$. 
\par
In terms of étale cohomology the functional equation above is a
consequence of the 
Poincaré duality (in this case, it also amounts to the Weil pairing for
the jacobian variety) which states that there is a natural non-degenerate
alternating pairing (``cup-product'')
\begin{equation}\label{eq-pairing}
H^1(\barre{C},\Zz_{\ell})\otimes H^1(\barre{C},\Zz_{\ell})\ra
\Zz_{\ell}(-1). 
\end{equation}
Note that this implies that the ``global'' geometric Frobenius $\frob$
of $\Fp_q$ acts on $H^1(\barre{C},\Zz_{\ell})$ as a symplectic
similitude for this pairing, with multiplicator $q$.
\par
\medskip
Here is now our first general result about the behavior of the splitting
fields in a suitable family, which significantly strengthens the
results of Chavdarov.

\begin{theorem}\label{th-chavdarov}
Fix an integer $g\geq 1$. Let $q=p^k$ and let $U/\Fp_q$ be a
geometrically irreducible smooth affine scheme of dimension $d\geq 1$
such that one of the following two conditions is satisfied:
\par
\quad \emph{(i)} $U$ is a curve, i.e, $d=1$,
\par
\quad \emph{(ii)} we have $p>2g+1$.
\par
Let $\pi\,:\,C\ra U$ be a proper smooth
family of projective curves of genus $g$ over $U$. Assume that for all
$\ell>L_0$ the geometric monodromy group of the integral sheaves 
$R^1\pi_!\Zz_{\ell}$ is the full symplectic group $Sp(2g)$. Then the
number  $N(U/\Fp_{q})$ of $u\in U(\Fp_q)$ such that the numerator 
\begin{equation}\label{eq-numzeta}
P_u=\det(1-T\frob\,\mid\, H^1(\barre{C}_u,\Qq_{\ell}))\in \Zz[T]
\end{equation}
of the zeta function of the curve $C_u=\pi^{-1}(u)$
is reducible or has splitting field with degree strictly less than
$2^gg!$ satisfies
$$
N(U/\Fp_q)\ll q^{d-\gamma}(\log q)
$$
for $\gamma=\frac{1}{4g^2+3g+5}$ in case~\emph{(i)} and 
$\gamma=\frac{1}{10g^2+6g+8}$ in case~\emph{(ii)}, where the implied
constant depends only on $L_0$, $g$ and $\barre{U}/\barre{\Fp}_q$.  
\end{theorem}

Here is another  almost equivalent way of phrasing this:
consider the zeta function 
$$
\tilde{Z}(s)=\exp\Bigl(\sum_{n\geq 1}{\frac{N(U/\Fp_{q^n})}{n}
q^{-ns}}\Bigr).
$$
It follows from Theorem~\ref{th-chavdarov}
that it extends to a holomorphic function
on the half-plane $\Reel(s)>d-\gamma$. 
\par
The second result is a uniform version (in terms of $g$) for the
families of hyperelliptic curves already introduced.

\begin{theorem}\label{th-chavdarov-2}
Let $g\geq 1$, $p\not=2$, $q=p^k$ with $k\geq 1$. Let $f\in\Fp_q[X]$
be a monic polynomial of 
degree $2g$ with distinct roots in $\barre{\Fp}_q$, $U\subset \Aa^1$
be the complement of the set of zeros of $f$ and denote by
$\pi\,:\,C\ra U$ the family of hyperelliptic curves of genus $g$ given by
$$
C_u\,:\,y^2=f(x)(x-u)
$$
completed by the section at $\infty$, with projection $\pi(x,y,u)=u$.
\par
Then the number $N(f,q)$ of $u\in U(\Fp_q)$ such that the polynomial
$$
P_u=\det(1-T\frob\,\mid\, H^1(\barre{C}_u,\Qq_{\ell}))\in \Zz[T]
$$
is either reducible or has splitting field with 
degree strictly smaller than $2^gg!$ satisfies
$$
N(f,q)\ll q^{1-\gamma}(\log q)
$$
for $\gamma=\frac{1}{4g^2+3g+5}$, where the implied constant is
absolute.
\end{theorem}

Note that the uniform bound in this last result is only non-trivial if
$g^2$ is somewhat smaller than $\log q$, precisely if $4g^2=(\log
q)e^{-f(q)}$ with $\log q=o(f(q))$.
Still, it is an uncommon feature to be able to say \emph{anything} for
this kind of problems in a situation where $g$ and $q$ increase
together, instead of having first $q\ra +\infty$ (compare with the
discussion in~\cite[Introduction]{katz-sarnak}). 
\par
\medskip
Since Theorem~\ref{th-chavdarov-2} is more delicate, we will start by
proving it in the  Section~\ref{sec-chav-unif} after some
common preliminaries; we will then quickly deal with the somewhat simpler
case of Theorem~\ref{th-chavdarov}. Here we just make a few additional
remarks which are of independent interest.
\par
First of all, since the estimate of Theorem~\ref{th-chavdarov-2} is (in 
particular) uniform in $q$, it can also be used in ``horizontal''
direction, i.e., with $q=p$ 
varying. For instance, we deduce the following quite easily:
\begin{proposition}\label{pr-sample1}
Let $g\geq 1$ be an integer, $f\in\Qq[X]$ be a polynomial of
degree $2g$ with distinct complex roots. For $n\in \Zz$ not a root of
$f$, let $C_n/\Qq$ be the hyperelliptic curve of genus $g$ with equation
$$
C_n\,:\, y^2=f(x)(x-n)
$$
and let $J_n$ be its Jacobian.
Then for $N\geq 3$, the number $S(N)$ of integers $n$ with
$|n|\leq N$ such that $J_n/\Qq$ is not simple up to isogeny satisfies 
$$
S(N)\ll N^{1/2-\delta}(\log\log N)
$$
where $\delta=\demi\tfrac{1}{4g^2+3g+5}$. The implied constant depends
on $g$ and the splitting field of $f$.
\end{proposition}

This should be compared with the individual global results
of~\cite[\S 6]{chavdarov}; our result is on average, but note that we
do not require any information on the image of the Galois
representations associated to $J_n$, and in particular we get results
valid for \emph{all} genus, independently of the endomorphism ring of
$J_n$ or any other global property.

\begin{proof}
Denote first by $Q_f$ the set of primes $p$
totally split in the splitting field of $f$. Notice that for any
$X\geq 2$ we have the sieving estimate 
$$
S(N)\leq |\{n\in\Zz\,\mid\, |n|\leq N\text{ and } n\mods{p}\notin
\Omega(p)\text{ for }p\in Q_f,\ p\leq X\}|
$$
where
$$
\Omega(p)=\{t\in\Zz/p\Zz\,\mid\,  f(t)\not\equiv 0\mods{p}\text{ and }
\det(1-T\frob_t\,\mid\,
H^1(\barre{C}_t,\Qq_{\ell}))\text{ is irreducible.}\}
$$
By Theorem~\ref{th-chavdarov-2} there exists a constant $C\geq 0$ such
that for all $p$ we have
\begin{equation}\label{eq-vertical-sieve}
|\Omega(p)|\geq p-Cp^{1-\gamma}(\log p)
\end{equation}
with $\gamma=\tfrac{1}{4g^2+3g+2}$. By the usual (strong) form of the
large sieve (see e.g.~\cite[Th. 6]{bombieri},~\cite[Th. 7.14]{ant}),
we have 
\begin{equation}\label{eq-ls-2}
S(N)\ll (N+X^2)J^{-1}
\end{equation}
where 
$$
J=\sumb_{q\leq X}{
\prod_{\stacksum{p\mid q}{p\in Q_f}}{\frac{|\Omega(p)|}{p-|\Omega(p)|}}},
$$
the $\flat$ sign indicating a sum restricted to squarefree numbers.
\par
Just taking primes into account we get by~(\ref{eq-vertical-sieve})
$$
J\gg \sum_{\stacksum{p\leq X}{p\in Q_f}}{
p^{\gamma}(\log p)^{-1}}
\gg X^{1+\gamma}(\log\log X)^{-1},
$$
by the Chebotarev density theorem, and taking $X=N^{1/2}$ the
proposition follows from~(\ref{eq-ls-2}). 
\end{proof}

Another corollary of the large sieve estimates is to families of
abelian varieties.

\begin{corollary}\label{cor-univ}
Let $q=p^k$ and $g\geq 1$ such that $p>2g+1$. Then the number
$N(g,q)$ of isomorphism classes of principally polarized abelian
varieties $A/\Fp_q$ such that the polynomial $\det(1-\frob T\,\mid\,
H^1(\barre{A},\Qq_{\ell}))$ is either reducible or has splitting field with
Galois group strictly smaller than $W_{2g}$ satisfies
$$
N(g,q)\ll q^{g(g+1)/2-\gamma}(\log q)
$$
where $\gamma=\frac{1}{10g^2+6g+8}$ and the implied constant depends
only on $p$ and $g$.
\end{corollary}

We will prove this at the same time as Theorem~\ref{th-chavdarov}.
\par
Finally, here is an (amusing, but not \emph{too} far-fetched)
illustration of what Theorem~\ref{th-chavdarov} gives which can not be
derived from~\cite{chavdarov}. 

\begin{proposition}\label{pr-sample2}
Let $g\geq 1$, $q=p^k$, $U/\Fp_q$ a non-empty open
subscheme of $\mathbf{G}_m/\Fp_q$, $\pi\,:\,C\ra U$ a family of smooth
projective curves of genus $g$ such that the geometric monodromy group
of  $R^1\pi_!\Zz_{\ell}$ is equal to $Sp(2g)$ for almost
all $\ell$.  Then for all $n$ large enough,  
there exists a \emph{primitive root} $t\in U(\Fp_{q^n})\subset
\Fp^{\times}_{q^n}$ such that $\det(1-T\frob_t\,\mid\,
H^1(\barre{C}_t,\Zz_{\ell}))$ has maximal Galois group.
\end{proposition}

\begin{proof}
Indeed, the set of $t\in \mathbf{G}_m(\Fp_{q^n})$ which are primitive
roots has cardinality $\varphi(q^n-1)$ and 
$$
\varphi(q^n-1)\gg \frac{q^n-1}{\log\log
  (q^n-1)}\gg \frac{q^n}{\log n+\log\log q}
$$
for $n\geq 1$ with an absolute implied constant (see
e.g.~\cite[Th. 328]{hardy-wright} for this 
standard estimate), and this lower bound is larger than the upper bound given
by Theorem~\ref{th-chavdarov} for those $t$ for which the numerator of
the zeta function of $C_t$ has small Galois group, if $n$ is large
enough. (So in fact, most primitive roots $t$ will have the desired property).
\end{proof}

\begin{remark}
For any $k$ coprime with $q$ we can find $n$ such that
$q^n-1\equiv 0\mods{k}$ and then
$$
\frac{\varphi(q^n-1)}{q^n-1}\leq \frac{\varphi(k)}{k}.
$$
Choosing suitable values of $k$, we see that the density of primitive
roots in $\Fp_{q^n}^{\times}$ is not bounded from below by any positive
constant. This means that Proposition~\ref{pr-sample2} can not be
proved without a quantitative form of Chavdarov's theorem.
\end{remark}

\section{Preliminaries for the proof of Chavdarov's theorem}
\label{sec-prelim-chav}

We start with some preliminaries related to the group $W_{2g}$ and to
setting up a sieve for characteristic polynomials of symplectic
similitudes. 
\par
From the description of $W_{2g}$ we see that there is an exact sequence
$$
1 \ra \{\pm 1\}^g\ra W_{2g}\fleche{p} \mathfrak{S}_g\ra 1
$$
where the second map just look at the permutation of the pairs, and
the kernel corresponds to just switching the elements of the pairs
without moving them. 
We also denote by $i$ the natural inclusion $i\,:\,W_{2g}\ra
\mathfrak{S}_{2g}$. 
\par
Our first lemma describes various ways of ensuring that a subgroup of
$W_{2g}$ is equal to $W_{2g}$.

\begin{lemma}\label{lm-separate}
Let $g\geq 1$ and $W\subset W_{2g}$ be a subgroup of $W_{2g}$. Assume
that one of the following conditions is true, where $i\,:\, W_{2g}\ra
\mathfrak{S}_{2g}$ is the embedding above:
\par
\quad\emph{(i)} For any conjugacy class $c\subset W_{2g}$, we have
$c\cap W\not=\emptyset$.
\par
\quad\emph{(ii)} The subgroup $i(W)$ contains a $2$-cycle, a
$4$-cycle, a $(2g-2)$-cycle and a $2g$-cycle.
\par
\quad\emph{(iii)} The subgroup $i(W)$ contains a transposition and
acts transitively on $\{1,\ldots, 2g\}$; moreover, the projection $p(W)$
contains a transposition and an $m$-cycle with $m>g/2$ prime.
\par
Then in all cases we have $W=W_{2g}$.
\end{lemma}

\begin{proof}
Case (i) is a standard result in finite group theory (see
e.g.~\cite[Lemma 5.8]{chavdarov}), which is in no way specific to
$W_{2g}$.
\par
Case (ii) is Lemma~2 of~\cite{dds}.
\par
For case (iii), observe first that the first condition already implies that
$W=W_{2g}$ if $g=1$. Otherwise we see that $p(W)$ acts
transitively on $\{1,\ldots, g\}$ and so with the second and third
conditions, we get $p(W)=\mathfrak{S}_g$ by
the result of Bauer given in~\cite[Lemma, p. 98]{gallagher}. Since
$i(W)$ contains a transposition, we deduce that $W=W_{2g}$ by
Lemma~5.5 of~\cite{chavdarov}.
\end{proof}

For Theorem~\ref{th-chavdarov}, we can use Case~(i) or
Case~(ii) of the lemma, but Theorem~\ref{th-chavdarov-2} requires the
finer Case~(iii), the point being that the conditions involve
``large'' subsets of $W_{2g}$. 
It seems to be an intriguing problem in combinatorial group theory to
determine how optimal this statement is.
In terms of $\mathfrak{S}_g$,
this means the following optimization problem: let $\Omega_1$, \ldots,
$\Omega_k$ be conjugacy-invariant subsets of $\mathfrak{S}_g$, and let
$$
\delta(\Omega_i)=\min \frac{|\Omega_i|}{g!}.
$$
How large can $\delta(\Omega_i)$ be if $(\Omega_i)$ are chosen so that
no proper subgroup of $\mathfrak{S}_g$ can intersect each $\Omega_i$?
Bauer's lemma gives two subsets with $\delta(\Omega_1,\Omega_2)\gg
1/\sqrt{g}$ (see~(\ref{eq-c3c4}) below).
\par
\medskip
To set up our sieve, it will be convenient to say that a polynomial
$f\in A[T]$ ($A$ any commutative ring) of degree $2g$ such that
$f(0)=1$ and
$$
T^{2g}f\Bigl(\frac{q}{T}\Bigr)=f(T),
$$
is \emph{$q$-symplectic of degree $2g$} ($q$ will
often be fixed and clear from the context, as will $g$).\footnote{\
The terminology ``self-reciprocal'' is often used when $q=1$.} Hence
the numerator of the zeta function of a curve $C/\Fp_q$ is
$q$-symplectic.
\par
We now prove a general result comparing a sieve related to
characteristic polynomials of elements  with multiplicator $q$ in the
finite group $SSp(2g,\Fp_{\ell})$ of symplectic similitudes
to the ``same'' sieve applied to all $q$-symplectic 
polynomials of degree $2g$. 
\par
Recall that we denote by $m(g)$ the multiplicator for a symplectic
similitude, i.e.,
$$
\langle gv,gw\rangle = m(g)\langle v,w\rangle.
$$

\begin{lemma}\label{lm-estim-general}
Let $g\geq 1$ and $\ell$ a prime. Put
$$
\Upsilon_{g,\ell}=\{f\in \Fp_{\ell}[T]\,\mid\, f\text{ is
  $q$-symplectic of degree $2g$}\}.
$$
\par
Let $\tilde{\Omega}(\ell)\subset \Upsilon_{g,\ell}$ be an arbitrary
subset of cardinality $\tilde{\omega}(\ell)$ and 
$$
\Omega(\ell)=\{g\in SSp(2g,\Fp_{\ell})\,\mid\,
m(g)=q,\text{ and }
\deg(1-Tg)\in\tilde{\Omega}(\ell)\},
$$
with cardinality $\omega(\ell)$.
\par
Then we have
$$
\omega(\ell)|Sp(2g,\Fp_{\ell})|^{-1}
\geq \tilde{\omega}(\ell)(\ell+1)^{-g}.
$$
\end{lemma}

\begin{proof}
We have
$$
\omega(\ell)=\sum_{f\in \tilde{\Omega}(\ell)}{
|\{g\in SSp(2g,\Fp_{\ell})\,\mid\, m(g)=q\text{ and }
\det(1-Tg)=f\}|}.
$$
The inner quantity, for given $f$, is exactly what is estimated by
Chavdarov in~\cite[Th. 3.5]{chavdarov}, in the proof of which it is
called $\Delta$. Using the formula at the bottom of page 159 of
loc. cit., we get
$$
\omega(\ell)=\frac{1}{\ell^g}\sum_{f\in \tilde{\Omega}(\ell)}{
|Sp(2g,\Fp_{\ell})|\frac{\ell^{\delta(f)}}{|C(A_f)(\Fp_{\ell})|}},
$$
where $A_f$ is a fixed semisimple element in $SSp(2g,\Fp_{\ell})$ with
multiplicator $q$ and characteristic polynomial $f$ (its existence
being proved in~\cite[Lemma 3.4]{chavdarov}), and $C(A_f)$ is the
centralizer of $A_f$ in $Sp(2g,\Fp_{\ell})$, $\delta(f)\leq g$ being its
dimension. Thus 
$$
\omega(\ell)|Sp(2g,\Fp_{\ell})|^{-1}=
\ell^{-g}\sum_{f\in \tilde{\Omega}(\ell)}{
\frac{\ell^{\delta(f)}}{|C(A_f)(\Fp_{\ell})|}}.
$$
By the formula of Nori at the top of page 160 of loc. cit., which
holds essentially because $C(A_f)$ is known to be a geometrically
irreducible variety of dimension $\delta(f)\leq g$, we have
$$
\omega(\ell)|Sp(2g,\Fp_{\ell})|^{-1}\geq 
\ell^{-g}\sum_{f\in \tilde{\Omega}(\ell)}{
\Bigl(1-\frac{1}{\ell+1}\Bigr)^{\delta(f)}}\geq
\ell^{-g}\Bigl(1-\frac{1}{\ell+1}\Bigr)^g\tilde{\omega}(\ell),
$$
as required.
\end{proof}

The next results are technical estimates which are only required in
this precise form for the proof of the uniform version of Chavdarov's
theorem. Easier versions (found in~\cite{gallagher},~\cite{dds})
suffice for Theorem~\ref{th-chavdarov}.
\par
Recall the following terminology: if $f$ is a monic polynomial of
degree $g$ in $\Zz[T]$ 
which factorizes modulo a prime $\ell$ as
$$
f=f_1\cdots f_r
$$
with the  $f_i$ coprime, irreducible, of degree $d_i\geq 1$, then one
says that the \emph{cycle type} (or the conjugacy class) associated to
$f$ is the conjugacy class in $\mathfrak{S}_g$ of elements which are
product of disjoint cycles of lengths $d_1$, \ldots, $d_r$.

\begin{lemma}\label{lm-omeg}
\emph{(i)} Let $g\geq 1$ and let $c$ be a conjugacy class in
$\mathfrak{S}_g$. For $\ell$ prime, let
$$
\hat{\Omega}_c(\ell)=
\{
f\in\Fp_{\ell}[T]\,\mid\, f\text{ is monic of degree $g$ and the cycle
type  associated to $f$ is $c$}\},
$$
and $\hat{\omega}_c(\ell)=|\hat{\Omega}_c(\ell)|$. Then we have
for $\ell>4g^2$
$$
\hat{\omega}_c(\ell)\geq \frac{|c|}{|\mathfrak{S}_g|}(\ell-1)^g
\Bigl(1-\frac{1}{\sqrt{\ell}}\Bigr)^g.
$$
\par
\emph{(ii)} Let $g\geq 1$ and for $\ell$ prime let $\omega_1(\ell)$ be
the number of $q$-symplectic irreducible polynomials in
$\Fp_{\ell}[T]$ of degree $2g$. 
Then for $\ell>4g^2$ we have
$$
\omega_{1}(\ell)\geq
\frac{\ell^g}{2g}\Bigl(1-\frac{1}{\ell}\Bigr)^g
-\ell^{g/2}.
$$
\par
\emph{(iii)} Let $g\geq 1$ and for $\ell$ prime let $\omega_2(\ell)$
be the number of $q$-symplectic polynomials of degree $2g$ 
which factorize as a product of an irreducible quadratic polynomial
and a product of irreducible polynomials of odd degrees. Then we have
for $\ell>4g^2$
$$
\omega_2(\ell)\geq \frac{\ell^g}{4g}\Bigl(1-\frac{1}{\ell}\Bigr)^g.
$$
\end{lemma}

\begin{proof}
We start with (i). If the conjugacy class $c$ is that consisting of
permutations with $r_i$ distinct $i$-cycles in their decomposition,
with $1\cdot r_1+\cdots+ g\cdot r_g=g$, then we have
$$
\frac{|c|}{|\mathfrak{S}_g|}=\prod_{1\leq i\leq g}
{\frac{1}{i^{r_i}r_i!}},
\text{ and }
\hat{\omega}_c(\ell)=\prod_{1\leq i\leq g}
{\binom{p(i,\ell)}{r_i}},
$$
where $p(i,\ell)$ is the number of irreducible monic polynomials of
degree $i$ in $\Fp_{\ell}[T]$. Now we claim that we have for all
$\ell>4g^2$ and $1\leq i\leq g$ the lower bounds
\begin{align}
p(1,\ell)&\geq \ell\Bigl(1-\frac{1}{\sqrt{\ell}}\Bigr)
\Bigl(1-\frac{1}{\ell}\Bigr)+g-1,
\label{eq-pil1}
\\
p(i,\ell)&\geq \frac{\ell^i}{i}\Bigl(1-\frac{1}{\ell}\Bigr)
+\frac{g}{i}-1\text{ for } 2\leq i\leq g.\label{eq-pil}
\end{align}
From this, which we prove below, we
derive for all $i$, $2\leq i\leq g$ and $r_i\leq g/i$ that
\begin{align*}
\binom{p(i,\ell)}{r_i}&=
\frac{p(i,\ell)(p(i,\ell)-1)\cdots (p(i,\ell)-r_i+1)}{r_i!}\geq
\frac{(p(i,\ell)-g/i+1)^{r_i}}{r_i!}\\ 
&\geq \Bigl(1-\frac{1}{\ell}\Bigr)^{r_i}\frac{1}{i^{r_i}r_i!}
\ell^{ir_i}
\end{align*}
and for $i=1$, $r_i\leq g$ that
\begin{align*}
\binom{p(1,\ell)}{r_1}&=
\frac{p(1,\ell)(p(1,\ell)-1)\cdots (p(1,\ell)-r_1+1)}{r_1!}
\geq \frac{(p(1,\ell)-g+1)^{r_1}}{r_1!}\\
&\geq \Bigl(1-\frac{1}{\sqrt{\ell}}\Bigr)^{r_1}
\Bigl(1-\frac{1}{\ell}\Bigr)^{r_1}\frac{1}{1^{r_1}r_1!}
\ell^{r_1}.
\end{align*}
Hence, putting these together, we get
\begin{align*}
\hat{\omega}_c(\ell)&\geq  
\Bigl(1-\frac{1}{\sqrt{\ell}}\Bigr)^{r_1}
\Bigl(1-\frac{1}{\ell}\Bigr)^{\sum r_i}\Bigl(\prod_{1\leq i\leq g}
{\frac{1}{i^{r_i}r_i!}}\Bigr)
\ell^{\sum ir_i}\\
&=\frac{|c|}{|\mathfrak{S}_g|}\Bigl(1-\frac{1}{\sqrt{\ell}}\Bigr)^{r_1}
\Bigl(1-\frac{1}{\ell}\Bigr)^{\sum r_i}\ell^g
\\
&\geq \frac{|c|}{|\mathfrak{S}_g|}\ell^g
\Bigl(1-\frac{1}{\ell}\Bigr)^{g}
\Bigl(1-\frac{1}{\sqrt{\ell}}\Bigr)^{g}
\end{align*}
as desired.
\par
Now we prove~(\ref{eq-pil1}) and~(\ref{eq-pil}). We use the well-known
formula of Dedekind
$$
p(i,\ell)=\frac{1}{i}\sum_{d\mid i}{\mu(d)\ell^{i/d}}.
$$
In particular
$$
p(1,\ell)=\ell-1\geq \ell\Bigl(1-\frac{1}{\sqrt{\ell}}\Bigr)
\Bigl(1-\frac{1}{\ell}\Bigr)+g-1
$$
for $\ell>4g^2$ by inspection. Similarly for $i=2$ we have
$$
p(2,\ell)=\frac{1}{2}(\ell^2-1)\geq \frac{1}{2}(\ell-1)^2+\frac{g}{2}-1
$$
if $\ell\geq g$. For $i\geq 3$ we use the lower bound
$$
p(i,\ell)\geq \frac{\ell^i}{i}-\ell^{i/2}
$$
(see e.g.~\cite[Lemma 3.1]{chavdarov}), so that it suffices to show
that
$$
\frac{\ell^i}{i}-\ell^{i/2}>\frac{\ell^i}{i}\Bigl(1-\frac{1}{\ell}\Bigr)
+\frac{g}{i}
$$
for $\ell>4g^2$. This is equivalent with
$$
\frac{1}{\ell}X^2-iX-g>0\text{ where } X=\ell^{i/2},
$$
and the quadratic polynomial has largest root equal to
$$
\alpha=\frac{\ell}{2}(i+\sqrt{i^2+4g/\ell})<2\ell i
$$
for $i\geq 3$ and $g<\sqrt{\ell}/2$. Hence for $i\geq 3$, $\ell^{1/2}>2g$
we have trivially $X=\ell^{i/2}>2\ell g\geq 2\ell i>\alpha$, so the quadratic
polynomial must be $>0$ when evaluated at $X$, which
gives~(\ref{eq-pil}). 
\par
Coming to (ii) we have (compare~\cite[Lemma 3]{dds}) the lower bound
$$
\omega_{1}(\ell)\geq p(g,\ell)-\frac{1}{2}p(g,\ell)-\ell^{g/2}.
$$
This is 
because we can count irreducible polynomials of degree $g$ in
$\Fp_{\ell}[T]$, minus those for which $f=T^gh(qT+T^{-1})$ is
reducible; in this case, $f$ is of the form $ch(T)T^gh(qT^{-1})$ (for
some normalizing constant $c\not=0$) where $h$ is irreducible of
degree $g$ and \emph{not} $q$-symplectic, with both $h$ and
$T^gh(qT^{-1})$ yielding (with proper normalization factor, so they
are distinct up to scalars by virtue of $h$ not being $q$-symplectic)
the same reducible $f$. From the irreducible $h$, we exclude the
$q$-symplectic ones by the trivial bound $\ell^{g/2}$ for their
number, hence the inequality above. 
\par
Using~(\ref{eq-pil}) for $i=g$ we get
$$
\omega_{1}(\ell)\geq 
\frac{\ell^g}{2g}\Bigl(1-\frac{1}{\ell^g}\Bigr)-\ell^{g/2}.
$$
Finally for (iii) we consider separately the case where $g$ is even or
when $g$ is odd. For even $g$, the number $\omega_2(\ell)$ is larger
than that of $q$-symplectic polynomials $f$ of degree $2g$ of the form
$$
f=f_1h_1h_2
$$
where $f_1$ is an irreducible quadratic $q$-symplectic polynomial, and
$h_1$, $h_2$ are irreducible of odd degree $g-1$ with, up to a constant,
$h_1=T^gh_2(qT^{-1})$. Counting the possibilities  we get
by~(\ref{eq-pil}) that 
$$
\omega_2(\ell)\geq \frac{\ell}{2}\Bigl(1-\frac{1}{\ell}\Bigr)
p(g-1,\ell)\geq \frac{\ell}{2}\Bigl(1-\frac{1}{\ell}\Bigr)
\frac{\ell^{g-1}}{g-1}\Bigl(1-\frac{1}{\ell}\Bigr)^{g-1}=
\frac{\ell^g}{2(g-1)}\Bigl(1-\frac{1}{\ell}\Bigr)^{g}
$$
hence a stronger result than claimed in this first case.
\par
The case where $g$ is odd is similar with polynomials of the form
$f=f_1f_2f_3h_1h_2$ with $f_1$ quadratic irreducible and
$q$-symplectic as before, $f_2=1-\alpha T$ for some $\alpha\not=0$ and
$f_3=1-q\alpha^{-1}T$, and $h_1$, $h_2$ as in the even case but now
with odd degree $g-2$. One gets a denominator $4(g-2)\leq 4g$ this
time.  
\end{proof}

\begin{remark}
As a by-product of these estimates, applying Gallagher's method 
we can derive a 
uniform version of his estimate~\cite{gallagher} for the
number $E_n(N)$ of monic polynomials 
$f=\sum{a_iT^i}\in\Zz[T]$ of degree $n\geq 1$ with height
$\max|a_i|\leq N$ and Galois group 
strictly smaller than $\mathfrak{S}_n$, namely
\begin{equation}\label{eq-gallagher-uniform}
E_n(N)\ll n^2(2N+1)^{n-1/2}(\log N)
\end{equation}
with an absolute implied constant (note $(2N+1)^n$ is the number
of polynomials with height $\leq N$).
\par
This is much more impressive than our Theorem~\ref{th-chavdarov-2}
because the gain in the exponent (namely, $\demi$) is independent of
the degree $n$ so the bound is non-trivial for $n$ as large as
$(2N)^{1/4}/(\log N)$. 
\par
Similarly for reciprocal polynomials, as treated in~\cite{dds}, denoting
by $\mathcal{E}_m(N)$ the number of monic
polynomials in $f\in\Zz[T]$ of degree $2m$ with height $\leq N$ such
that $T^mf(T^{-1})=f$ we get
$$
\mathcal{E}_m(N)\ll m^2(2N+1)^{m-1/2}(\log N)
$$
with an absolute implied constant. 
\par
Note that in the two papers quoted, the
fundamental large-sieve inequality is not
uniform in $n$ (resp. $m$) as stated, but becomes so if one replaces it
by the form given in~\cite[Th. 1]{huxley}, with some obvious changes.
For instance in~\cite[eq. (5)]{gallagher}, the term
$(N^n+x^{2n})$ in the right-hand side must be replaced by
$(\sqrt{2N+1}+x)^{2n}$. To make this innocuous one may take
$x=n^{-1}\sqrt{2N+1}$ instead of $x=\sqrt{N}$, which leads to
$\pi(x)^{-1}\ll n(2N+1)^{-1/2}(\log N)$, hence the ``extra'' power of $n$
in~(\ref{eq-gallagher-uniform}) compared with the ``density'' of
irreducible polynomials modulo $\ell$ (i.e, about $n^{-1}$) which are
used to sieve the reducible ones.
\end{remark}

\section{Proof of the uniform version of Chavdarov's
theorem}\label{sec-chav-unif} 

We can now start the proof of Theorem~\ref{th-chavdarov-2} itself.
We will apply Proposition~\ref{pr-large-sieve} with the following
data: in addition to $U$, which is a smooth geometrically connected
affine curve over $\Fp_q$, we take the sheaves
$\sheafm{F}=R^1f_!\Fp_{\ell}$ for $\ell\in \Lambda$,
where $\Lambda$ is the set of odd primes 
\par
These sheaves are of course obtained by reduction modulo $\ell$ from 
the compatible system $\sheaf{F}_{\ell}=R^1f_!\Zz_{\ell}$. 
The existence of the symplectic
pairing~(\ref{eq-pairing}) implies that
the arithmetic monodromy group of $\sheafm{F}$ can be seen as a
subgroup of $SSp(2g,\Fp_{\ell})$, and for any $u\in U(\Fp_q)$,
the image of $\frob_u$ has multiplicator $q$. 
\par
The most crucial point is that for $\ell>2$, J-K. Yu (unpublished) has
shown that the geometric monodromy group for $\sheafm{F}$ is equal to
$Sp(2g,\Fp_{\ell})$. Then the sheaves $(\sheafm{F})$ are also \disj\ as a
consequence of Goursat's Lemma (see Corollary~\ref{cor-goursat}, (1)),
and by Lemma~\ref{lm-kappa} we have~(\ref{eq-kappa}) with $\kappa=2$.
And finally, Lemma~\ref{lm-sizes}, (3) gives us~(\ref{eq-s}) with
$s=2g^2+g+1$, $t=g+1$ and $c_1=1$, $c_2=6^{g}$.
\par
Thus all conditions needed to apply Proposition~\ref{pr-large-sieve}
(in the case of a one-parameter family) are valid, and it remains to
set up the sieving problem. The principle for this is exactly the same
as the one introduced by Gallagher for polynomials with integer
coefficients and bounded height~\cite{gallagher}.
\par
As in Lemma~\ref{lm-estim-general}, for any choice of sets
$\tilde{\Omega}(\ell)\subset \Upsilon_{g,\ell}$ defined for $\ell>2$
we let 
$$
\Omega(\ell)=\{g\in SSp(2g,\Fp_{\ell})\,\mid\,
m(g)=q,\text{ and }
\deg(1-Tg)\in\tilde{\Omega}(\ell)\}.
$$
Applying Proposition~\ref{pr-large-sieve} (see~(\ref{eq-sieve})) to
such a sieving problem, we have
\begin{equation}
\label{eq-apply-ls}
|S(U,\Omega;L)|\leq (2q+4gq^{1/2}(6L)^{A})P(L)^{-1},
\end{equation}
where $A=2g^2+3g/2+5/2$ and
\begin{equation}\label{eq-density1}
P(L)=\sum_{2<\ell\leq L}{\omega(\ell)|G^g_{\ell}|^{-1}},
\end{equation}
(which shouldn't be confused with the polynomials $P_u$),
and here we have taken the constant $C=4g$ by looking at
Proposition~\ref{pr-boundcurve}
and~(\ref{eq-delta-curve}),~(\ref{eq-delta-curve-2}) since
$1-\chi(\barre{U},\barre{\Qq}_{\ell})=2g$ and it is known
that all the sheaves $\sheaf{F}_{\ell}$ are tamely ramified
(by~\cite[Lemma 10.1.12]{katz-sarnak}) so that 
the contribution $w$ of the Swan conductors in~(\ref{eq-ww})
vanishes. Moreover, by Lemma~\ref{lm-estim-general} we have
\begin{equation}\label{eq-basic-lower-bound}
P(L)\geq \sum_{2<\ell\leq L}{\tilde{\omega}(\ell)(\ell+1)^{-g}}.
\end{equation}
\par
Now we must show how to use this sieve estimate to study the
characteristic polynomials $P_u$. For this we need to recall the
following two facts, the first of which is classical, while the second
is much deeper: 
\par
\quad (i) if $f\in\Zz[T]$ is a polynomial of degree $d$ that
factorizes in $\Fp_{\ell}[T]$ as a product 
of coprime polynomials $f_1\cdots f_r$, with $f_i$ irreducible of
degree $\deg f_i=d_i$, then the Galois group of $f$, seen as a permutation
group of the complex roots of $f$, contains a cycle $c$ of type $(d_1,\ldots,
d_r)$, i.e., a product of disjoint cycles of respective length $d_1$,
\ldots, $d_r$ (see e.g.~\cite[\S 61]{vdw},~\cite[p. 302]{jacobson}).
\par
\quad (ii) the reduction modulo a prime $\ell$ of a polynomial $P_u$
(the numerator of the zeta function of the curve $C_u=\pi^{-1}(u)$)
is the characteristic polynomial of $\frob_u$ acting on
$\sheafm{F}$ (see~\cite[Fonctions L mod. $\ell^n$]{sga4half}, or use
the fact that
$$
\det(1-T\frob_{u}\,\mid\, R^1\pi_!\Zz_{\ell})=
\det(1-T\frob \,\mid\, H^1(\barre{C}_u,\Zz_{\ell}))=P_u
$$
by~(\ref{eq-cohomo-interp}), and reduce modulo $\ell$). 
\par
Thus (ii) allows us to control the reduction of a polynomial $P_u$,
while (i) tells us that the reduction gives information on the Galois
group of $P_u$.
\par
In particular, for any sieving sets $\Omega(\ell)\subset
SSp(2g,\Fp_{\ell})$, an element  
$u\in S(U,\Omega;L)$ will have the property that the Galois group of
$P_u$, seen as a subgroup of $\mathfrak{S}_{2g}$, does not contain a
cycle $c$ associated to an $f\in \Omega(\ell)$, where $\ell$ ranges
over primes $2<\ell\leq L$. 
\par
If we have finitely many sieving sets $\Omega_{i}$, $1\leq
i\leq m$, defined by the condition that the cycle associated to
$\det(1-Tg)$ is in a certain set $c_i$ of conjugacy classes, and if
moreover those $c_i$ have the property that the only subgroup
$W\subset W_{2g}$ containing an element of each $c_i$ is $W_{2g}$, 
then it follows that the set of exceptional $u\in U(\Fp_q)$ with $P_u$
having small Galois group will be a subset of the union of the
$S(U,\Omega_{i};L)$. So in such a situation we have
\begin{equation}\label{eq-separate}
N(U/\Fp_q)\leq S(U,\Omega_{1};L)+\cdots+S(U,\Omega_{m};L)
\leq (2q+4gq^{1/2}L^{A})\sum_{1\leq i\leq m}{P_{i}(L)^{-1}}.
\end{equation}
\par
Lemma~\ref{lm-separate} describes three
possible choices of sets $c_i$; however, the first and the second 
involve some $c_i$ which are ``too small'', so the dependency on $g$
in the estimate for $P_c(L)$ is bad (they are perfectly suitable for
fixed $g$). Thus we use Case~(iii) of
Lemma~\ref{lm-separate}. Precisely, we have $m=4$ and the four sets
$\Omega_i$ can be described as follows:
\par
\quad (i) $\Omega_1$ is the set of irreducible polynomials $f\in
\Upsilon_{g,\ell}$. 
\par
\quad (ii) $\Omega_2$ is the set of polynomials $f\in \Upsilon_{g,\ell}$
which factorize as a product of an irreducible quadratic polynomial
and a product of irreducible polynomials of odd degrees.
\par
To define $\Omega_3$ and $\Omega_4$, we recall that any
$f\in\Upsilon_{g,\ell}$ can be written uniquely
$$
f=T^gh(qT+T^{-1})
$$
where $h\in\Fp_{\ell}[T]$ is a monic polynomial of degree $g$. 
\par
\quad (iii) $\Omega_3$ is the set of $f\in\Upsilon_{g,\ell}$ such that the
corresponding $h$ has an irreducible factor of prime degree $>g/2$.
\par
\quad (iv) $\Omega_4$ is the set of $f\in\Upsilon_{g,\ell}$ such that the
corresponding $h$ has a single quadratic irreducible factor and no
other irreducible factor of even degree.
\par
\medskip
We claim that those sets do allow us to sieve the exceptional elements
$u$. Indeed, spelling out again the
relation between the factorization of $P_u$ modulo $\ell$ and the
existence of elements in the Galois group of $P_u$ with the associated
cycle type, we see that:
\par
\quad (i) If $P_u$ is reducible then $u\in S(U,\Omega_{1};L)$.
\par
\quad (ii) If $P_u$ is irreducible but the Galois group $W$ does not
contain a transposition, then $u\in S(U,\Omega_{2};L)$, since having
$P_u\mods{\ell}\in \tilde{\Omega}_{2}(\ell)$ implies that $W$
contains an element with cycle type consisting of one $2$-cycle and
further cycles of odd length, a power of which will be a transposition.
\par
For the next two facts, notice that the cycle in $\mathfrak{S}_g$
associated to the polynomial $Q_u$ such that $P_u=T^gQ_u(qT+T^{-1})$
is the image by the map $p\,:\, W_{2g}\ra \mathfrak{S}_g$ of the cycle
associated to $P_u$.
\par
\quad (iii) If $P_u$ is irreducible but $p(W)$ does not contain a
cycle of prime order $m>g/2$, then $u\in S(U,\Omega_{3};L)$.
\par
\quad (iv)  If $P_u$ is irreducible but $p(W)$ does not contain a
transposition, then $u\in S(U,\Omega_{4};L)$ (as in Case~(ii)
previously). 
\par
By Case~(iii) of Lemma~\ref{lm-separate}, the $u\in U(\Fp_q)$ that we
wish to exclude are therefore in the union of the $S(U,\Omega_i;L)$,
and we conclude that
\begin{equation}\label{eq-borne-n}
N(U/\Fp_q)\leq S(U,\Omega_{1};L)+\cdots+S(U,\Omega_{4};L)
\leq 4(2q+4gq^{1/2}(6L)^{A})
\Bigl(\min\sum_{1\leq i\leq m}{P_{i}(L)}\Bigr)^{-1}.
\end{equation}
\par
It remains to give appropriate lower bounds of $P_i(L)$. For
$\Omega_{3}$ and $\Omega_{4}$, since the correspondence between
polynomials $f\in \Upsilon_{g,\ell}$ and the $h\in \Fp_{\ell}[T]$ such
that $f=T^gh(qT+T^{-1})$ is one-to-one, we can count the corresponding $h$
by Lemma~\ref{lm-omeg}, applied to the cycle types (i.e.,
conjugacy classes) in $\mathfrak{S}_g$ associated to the polynomials
in $\Omega_{i}$. For $\ell>4g^2$ and $i\in\{3,4\}$, denoting by $C_i$  
the set of elements in $\mathfrak{S}_g$ having the associated cycle
type, we get
$$
\tilde{\omega}_{i}(\ell)\geq \frac{|C_i|}{|\mathfrak{S}_g|}(\ell-1)^g
\Bigl(1-\frac{1}{\sqrt{\ell}}\Bigr)^g,
$$
and thus for $L>4g^2$ we have
$$
P_{i}(L)\geq \frac{|C_i|}{|\mathfrak{S}_g|}\sum_{4g^2<\ell\leq L}{
\Bigl(\frac{\ell-1}{\ell+1}\Bigr)^g
\Bigl(1-\frac{1}{\sqrt{\ell}}\Bigr)^g}.
$$
By the mean-value theorem we have for any $\ell\geq 2$
$$
\Bigl(\frac{\ell-1}{\ell+1}\Bigr)^g
\Bigl(1-\frac{1}{\sqrt{\ell}}\Bigr)^g
=1-gh(\ell)+O(g^2h(\ell)^2)
$$
with
$$
h(\ell)=\frac{2}{\ell+1}+\frac{1}{\sqrt{\ell}}-\frac{2}{\sqrt{\ell}(\ell+1)},
$$
and an absolute implied constant. Inserting this in the sum and using
the prime number theorem we get for $L>4g^2$ that
\begin{equation}\label{eq-lb-2}
P_{i}(L)\geq \frac{|C_i|}{|\mathfrak{S}_g|}
\{\pi(L)+O(g\sqrt{L}(\log L)^{-1}+g^2\log\log L)\},
\end{equation}
with an absolute implied constant.
\par
By~\cite[p. 99]{gallagher} (where our $C_3$ is denoted $P$ and $C_4$
is denoted $T$), we have for $g\geq 1$
\begin{equation}\label{eq-c3c4}
\frac{|C_3|}{|\mathfrak{S}_g|}\gg \frac{1}{\log 2g}\text{ and }
\frac{|C_4|}{|\mathfrak{S}_g|}\gg \frac{1}{\sqrt{g}}.
\end{equation}
Using~(\ref{eq-lb-2}), this gives the lower
bounds 
\begin{equation}\label{eq-bc3c4}
P_{3}(L)\gg \frac{1}{\log 2g}L(\log L)^{-1},\text{ and }
P_{4}(L)\gg \frac{1}{\sqrt{g}}L(\log L)^{-1}
\end{equation}
with absolute implied constants for $L\gg g^2(\log 2g)$
(i.e. for $L\geq \alpha_1 g^2(\log 2g)$, where the absolute constant
$\alpha_1$ can be specified from the implied constants
in~(\ref{eq-lb-2}) and~(\ref{eq-c3c4})).
\par
Coming to $\Omega_1$, we have by (ii) of Lemma~\ref{lm-omeg} that for
$\ell>4g^2$  
$$
\tilde{\omega}_{1}(\ell)\geq
\frac{\ell^g}{2g}\Bigl(1-\frac{1}{\ell^g}\Bigr)
-\ell^{g/2}
$$
so by~(\ref{eq-basic-lower-bound}),  the prime number theorem and the
mean-value theorem as before we get for $L>4g^2$ that
$$
P_{1}(L)\geq \frac{1}{2g}(\pi(L)+O(g\log\log L+g^2+\sqrt{L}))
$$
with an absolute implied constant, and hence for $L\gg g^2(\log 2g)$,
we have
\begin{equation}\label{eq-bc1}
P_{1}(L)\gg \frac{1}{g}L(\log L)^{-1}
\end{equation}
with absolute implied constant.
\par
Finally by (iii) of Lemma~\ref{lm-omeg} we have for $\ell>4g^2$
$$
\tilde{\omega}_{2}(\ell)\geq
\frac{1}{4g}\Bigl(1-\frac{1}{\ell}\Bigr)^g
\text{ and }
P_{2}(L)\geq \frac{1}{4g}(\pi(L)+O(g\log\log L+g^2))
$$
and for $L\gg g^2(\log 2g)$ we obtain also
\begin{equation}\label{eq-bc2}
P_{2}(L)\gg \frac{1}{g}L(\log L)^{-1}
\end{equation}
with absolute implied constant.
\par
Altogether from~(\ref{eq-borne-n}),~(\ref{eq-bc3c4}),~(\ref{eq-bc1})
and~(\ref{eq-bc2}) we get
$$
N(U/\Fp_{q})\ll g^2(2q+q^{1/2}(6L)^A)L^{-1}(\log L)
$$
with an absolute implied constant, which can in fact be
chosen so that 
the inequality is valid for all $L\geq 2$ and $g\geq 1$, since it becomes
trivial for $g^2\gg L(\log L)^{-1}$. Choosing
$6L=q^{(2A)^{-1}}=q^{(4g^2+3g+5)^{-1}}$, with $\log L\ll g^{-2}\log q$,
this gives the announced uniform estimate
$$
N(U/\Fp_{q})\ll q^{1-\gamma}(\log q)
$$
with $\gamma=(4g^2+3g+5)^{-1}$, and an absolute implied constant.

\section{Proof of the general version of Chavdarov' Theorem}

We will now quickly prove Theorem~\ref{th-chavdarov}, only
highlighting the points where the proof is different from that of 
the previous section. The first step is to check that we can always 
apply Proposition~\ref{pr-large-sieve} to the data consisting of
$U/\Fp_q$ and the family of sheaves $\sheafm{F}=R^1\pi_!\Fp_{\ell}$,
defined for a subset $\Lambda$ of primes $\ell>L_0$. 
\par
Since our
assumption is that for $\ell>L_0$ the geometric monodromy group of
$R^1\pi_!\Zz_{\ell}$ is the symplectic group $Sp(2g)$ (as algebraic
group over $\Qq_{\ell}$), we must show that this implies that the
monodromy group modulo $\ell$ is often large. (A priori, for fixed
$\ell$, the assumption only implies that the index of the image of
$\pi_1(\barre{U},\barre{\eta})$ in $Sp(2g,\Zz/\ell^{\nu}\Zz)$ is
bounded for $\nu\geq 1$, but does not say anything for $\nu=1$).
However, we can appeal to a result of Larsen~\cite[Th. 3.17]{larsen}
which implies that for a set of primes $\Lambda_0$ of (natural)
density $1$, we do have $G_{\ell}^g=Sp(2g,\Fp_{\ell})$ because the
sheaves come from a compatible system. (Precisely, in the notation of
loc. cit., apply Th. 3.17 with
$\mathcal{G}=\pi_1(\barre{U},\barre{\eta})$, $\rho_{\ell}$ 
corresponding to $R^1\pi_!\Zz_{\ell}$, so that by assumption
$G_{\ell}=Sp(2g)/\Qq_{\ell}$, which is connected and simply connected
so $G_{\ell}^{sc}=Sp(2g)$, and look at the first few lines of the
proof of Th. 3.17 to make sure that the statement there involving
``hyperspecial maximal compact subgroups'' does imply that the
geometric monodromy group of the reduction of $\rho_{\ell}$ is
$Sp(2g,\Fp_{\ell})$ for $\ell$ in a set of density $1$; note also that
Larsen's result is quite deep as it depends on the classification of
simple finite groups). 
\par
Then, as before, Corollary~\ref{cor-goursat} implies that the sheaves
$\sheafm{F}$ for $\ell\in\Lambda_0$ are
\disj\ in all cases (if $\ell\geq 5$) by the assumption on the
geometric monodromy groups, and~(\ref{eq-kappa}) holds with $\kappa=2$
by Lemma~\ref{lm-kappa}, (2).
\par
Since $\sheafm{F}$ is obtained by reduction of the
compatible system 
$\sheaf{F}_{\ell}=R^1\pi_!\Zz_{\ell}$, case (i) of
Proposition~\ref{pr-large-sieve} is applicable if $U$ is a curve,
with $\Lambda$ consisting of all the primes in $\Lambda_0$.
\par
If $U$ is not a curve but $p>2g+1$, we use the following simple lemma
(compare~\cite[Lemma 7.5.1]{katz-tw}) :

\begin{lemma}\label{lm-gl}
Let $r\geq 1$. For any $p>r+1$, there
exists $\alpha\in(\Zz/p\Zz)^{\times}$ such that 
the order of $GL(r,\Fp_{\ell})$ is prime to $p$ for any prime
$\ell\equiv\alpha\mods{p}$.
\end{lemma}

\begin{proof}
The order of $GL(r,\Fp_{\ell})$ is
$$
|GL(r,\Fp_{\ell})|=\ell^{r(r-1)/2}\prod_{1\leq i\leq r}{(\ell^{i}-1)}
$$
so the condition will hold whenever the order of $\alpha$ modulo $p$
is $>r$. If $p>r+1$, a primitive root modulo $p$ will certainly work.
\end{proof}

For $p>2g+1$, we can apply the second case of Theorem~\ref{th-main}
and Proposition~\ref{pr-large-sieve} with $\Lambda$ consisting of
primes in $\Lambda_0$ which are congruent modulo $p$ to the  $\alpha$
given by this lemma for $r=2g$. This set has positive density among the
primes because $\Lambda_0$ has density $1$.
\par
We can now define sieving sets analogous to the previous $\Omega_i$,
$1\leq i\leq 4$, for 
$U/\Fp_{q}$ and we have~(\ref{eq-borne-n}). Since we consider $g$ to
be fixed here, we can rewrite~(\ref{eq-basic-lower-bound}) as
$$
P_{i}(L)\gg \sum_{\stacksum{L_0<\ell\leq L}{\ell\in\Lambda}}
{\tilde{\omega}_{i}(\ell)\ell^{-g}}
$$
where the implied constant depends on $g$. Then we obtain
$$
P_{i}(L)\gg \pi(L)
$$
for $L>L_0$, the implied constant depending on $g$ and $p$ in Case
(ii) (through the density of $\Lambda$), either from
Lemma~\ref{lm-omeg} or directly from Dedekind's formula used in its
proof. 
\par
Hence we obtain by Proposition~\ref{pr-large-sieve}
$$
N(U/\Fp_q)\ll (q^d+Cq^{d-1/2}L^A)L^{-1}(\log L),
$$
for $L>L_0$, where the implied constant depends on $g$, and on $p$ in Case
(ii). If we take
$$
L=q^{(4g^2+3g+5)^{-1}},\text{ in case (i)},\quad\quad
L=q^{(10g^2+6g+8)^{-1}},\text{ in case (ii)}
$$
we have $q^{d-1/2}L^A=q^{d}$, hence
$$
N(U/\Fp_q)\ll q^{d}L^{-1}(\log L)\ll q^{d-\gamma}(\log q)
$$
for $\gamma=\frac{1}{4g^2+3g+5}$ (resp
$\gamma=\frac{1}{10g^2+6g+8}$), as desired. 
\par
\medskip
The proof of Corollary~\ref{cor-univ} is similar, applying first the
large sieve to 
a suitably rigidified moduli space $A_{g}/\Fp_q$ of
principally polarized abelian varieties over $\Fp_q$, for instance the
moduli space $A_{g,3\mathcal{L}}$
of~\cite[11.3]{katz-sarnak}. Strictly speaking 
we need to restrict to a smooth connected affine subscheme $U\subset
A_{g,3\mathcal{L}}$, but this is not a problem as observed in the remarks
after Theorem~\ref{th-main} (see Remark~4). Over $U$ we have a
universal family $\pi\,:\,\mathcal{A}_{g,3\mathcal{L}}\ra U$ 
and we take the sheaves $\sheaf{F}_{\ell}=R^1\pi_!\Zz_{\ell}$ and
their reductions $\sheaf{F}_{\ell}/\ell\sheaf{F}_{\ell}$. The
monodromy groups are as large as possible because already this is the
case for the families of (canonically principally polarized) jacobians
of the hyperelliptic curves of genus $g$ of
Theorem~\ref{th-chavdarov-2}. After applying
Proposition~\ref{pr-large-sieve} to the same 
sieving problem as in Theorem~\ref{th-chavdarov}, we go (with the same
saving) from
the number of ``exceptional'' principally polarized abelian varieties
with a $3\mathcal{L}$-structure to the number $N(g,q)$ by dividing out
by the free rigidifying parameters and considering the situations with
extra automorphisms, as done in~\cite[11.3]{katz-sarnak} for 
the case of curves.

\end{document}